\documentclass[11pt]{article}
\usepackage{amssymb,latexsym,amsmath,amsbsy,amsthm,amsxtra,amsgen,dsfont,pdfsync,graphicx,bm,mathtools,tikz,longtable, tikz-cd}
\usepackage{extpfeil}
\usepackage{mathrsfs}
\usepackage{bbm}
\usepackage{indentfirst}
\usepackage{cite}
\usepackage{mathrsfs} 
\usepackage{xcolor}
\usepackage{url}
\usepackage{authblk}
\usepackage{blindtext}
\usepackage{booktabs}
\usepackage{titling}
\usepackage{comment}
\usepackage[shortlabels]{enumitem}
\usepackage[all]{xy}
\usepackage{scalerel,stackengine}

\stackMath
\newcommand\reallywidehat[1]{%
\savestack{\tmpbox}{\stretchto{%
  \scaleto{%
    \scalerel*[\widthof{\ensuremath{#1}}]{\kern-.6pt\bigwedge\kern-.6pt}%
    {\rule[-\textheight/2]{1ex}{\textheight}}
  }{\textheight}%
}{0.5ex}}%
\stackon[1pt]{#1}{\tmpbox}%
}
\predate{}
\postdate{}

\usetikzlibrary{decorations.pathreplacing}

\newfont{\msbmsm}{msbm10 at 8pt}

\newtheorem{Theo}{Theorem}[section]
\newtheorem{Lemma}[Theo]{Lemma}
\newtheorem{Cor}[Theo]{Corollary}
\newtheorem{Prop}[Theo]{Proposition}

\newtheorem{Dfn}[Theo]{Definition}

\newtheorem{Notation}[Theo]{Notations}

\newcommand{\qedwhite}{\hfill \ensuremath{\Box}}

\title{\textbf{Adic Sheafiness of} $\mathbf{A}_{\textnormal{inf}}$ \textbf{Witt Vectors over Perfectoid Rings}}
\author{Zongze Liu \thanks{\it E-mail address: zol013@ucsd.edu}}
\affil{\fontsize{9}{10.8}\itshape{Department of Mathematics, University of California San Diego, La Jolla, CA 92093-0112, USA}}

\date{}

\begin{document}
\maketitle

\begin{abstract}
For $(R, R^{+})$ an analytic perfectoid ring in char $p$, let $\mathbf{A}_{\textnormal{inf}}(R^{+})$ be the ring of Witt vectors with the induced topology from $(R, R^{+})$. We prove that $\textnormal{Spa}(\mathbf{A}_{\textnormal{inf}}(R^{+}),\mathbf{A}_{\textnormal{inf}}(R^{+}))$ is a sheafy adic space and its structure sheaf is acyclic. We first show $\mathbf{A}_{\textnormal{{inf}}}(R^{+})$ is a stably uniform Banach ring. The `stably uniform implies sheafy' argument is applied to Tate Huber rings in \cite{bv} and is generalized to analytic Huber rings in \cite{aws}. Here we show that the `stably uniform implies sheafy' argument in \cite{aws} can be applied to general stably uniform Banach rings whose underlying topological ring is a Huber ring. Finally we show the equivalence of categories of vector bundles over $\textnormal{Spa}(\mathbf{A}_{\textnormal{inf}}(R^{+}),\mathbf{A}_{\textnormal{inf}}(R^{+}))$ and finite projective modules over $\mathbf{A}_{\textnormal{inf}}(R^{+})$.
\end{abstract}

\section{Introduction}
Let $(R, R^{+})$ be a char $p$ perfectoid ring, $\mathbf{A}_{\textnormal{inf}}(R^{+})$ is the ring of Witt vectors $W(R^{+})$ equipped with the topology induced from $(R, R^{+})$ as defined in Notations \ref{ainf topology}. 
$\mathbf{A}_{\textnormal{inf}}$ plays a pivotal role in $p$-adic Hodge theory and $p$-adic geometry. $\mathbf{A}_{\textnormal{inf}}$ is first defined for a complete nonarchimedean perfect field $K$ equipped with a nontrivial valuation $v$ by Fontaine to construct the $p$-adic period rings for establishing the comparison isomorphisms in $p$-adic Hodge Theory. Later Bhatt-Morrow-Scholze\cite{bms} defined an $\mathbf{A}_{\textnormal{inf}}$-valued cohomology theory to reinterpret and refine the crystalline comparison isomorphism. In the theory of perfectoid spaces, $\mathbf{A}_{\textnormal{inf}}$ is instrumental for establishing the perfectoid correspondence between characteristic 0 and characteristic $p$ perfectoids. For a perfectoid ring $(R, R^{+})$ in characteristic $p$, $\mathbf{A}_{\textnormal{inf}}(R^{+})$ classifies the set of characteristic 0 untilt perfectoid rings $(S, S^{+})$ of $(R, R^{+})$ with its primitive elements: given a primitive element $z \in \mathbf{A}_{\textnormal{inf}}(R^{+})$, a untilt of $R$ is given by $(S, S^{+})$ where $S:= W^{b}(R)/(z)W^{b}(R) \ \textnormal{and} \ S^{+}:=\mathbf{A}_{\textnormal{inf}}(R^{+})/(z)$. Then the adic space $\textnormal{Spa}(S,S^{+})$ is canonically identified with  $\textnormal{Spa}(R,R^{+})$.The Fargues-Fontaine curve \cite{ff}, a central object in $p$-adic Hodge theory, was first defined as a schematic curve using algebraically closed nonarchimedean fields and then perfectoid fields. The relative and adic Fargues-Fontaine curve was introduced in \cite{kl} over Tate perfectoid rings and over analytic perfectoid rings in \cite{aws}. The relative Fargues-Fontaine curve over $\textnormal{Spa}(R,R^{+})$, $\mathbf{FF}_{R}$, is defined as the quotient space of an analytic open locus of $\textnormal{Spa}(\mathbf{A}_{\textnormal{inf}}(R^{+}),\mathbf{A}_{\textnormal{inf}}(R^{+}))$ by the action of the cyclic group generated by the Frobenius map $\varphi$. Moreover, a GAGA theorem between the analytic adic curve and the schematic curve was established in \cite{kl}. 
The relative Fargues-Fontaine curve $\mathbf{FF}_{R}$ is a sheafy adic space. When $R$ is the tilt of a complete algebraically-closed nonarchimedean field $\mathbf{C^{\flat}}$, the category of vector bundles over $\mathbf{FF}_{\mathbf{C^{\flat}}}$ with an Frobenius isomorphism is equivalent to the category of Breuil-Kisin-Fargues modules over $\mathbf{A}_{\textnormal{inf}}(R^{+})$ and the category of mixed characteristic shtukas with one leg, which is the starting point of the geometrization of $p$-adic local Langlands program, as summarized in \cite{berkeley}. \par 
The main result of our paper proves the sheafiness of $\textnormal{Spa}(\mathbf{A}_{\textnormal{inf}}(R^{+}),\mathbf{A}_{\textnormal{inf}}(R^{+}))$:
\begin{Theo}
Let $(R, R^{+})$ be an analytic perfectoid ring in char $p$. $\textnormal{Spa}(\mathbf{A}_{\textnormal{inf}}(R^{+}), \mathbf{A}_{\textnormal{inf}}(R^{+}))$ is a sheafy adic space and $\mathcal{O}_{\textnormal{Spa}(\mathbf{A}_{\textnormal{inf}}(R^{+}), \mathbf{A}_{\textnormal{inf}}(R^{+}))}$ is an acyclic sheaf.
\end{Theo}

The nonarchimedean analytic geometry of Witt vectors was first studied in \cite{K13} for p-typical Witt vector rings $W(S)$ equipped with the $p$-adic norm where $S$ is a perfect $\mathbb{F}_{p}$-algebra with the trivial norm. Working with Berkovich spaces, \cite{K13} already shows the homeomorphism of topological spaces that underlies the perfectoid correspondence. Despite the importance of $\mathbf{A}_{\textnormal{inf}}$ in the adic space geometry of perfectoid spaces and the relative Fargues-Fontaine curves, the adic space geometry of  $\mathbf{A}_{\textnormal{inf}}$ itself remains previously largely unexplored. This is due to the fact that the streamlined proofs of generalizations of the Tate acyclicity theorems to adic spaces as in \cite{bv} and \cite{aws} require the Huber rings to contain a pseudo-uniformizer or, more generally, be analytic. A Huber ring $(A, A^{+})$ is analytic if all of its valuations are analytic, i.e. the kernel of the valuation does not contain open ideals, or equivalently, the ideal of definition of $A$ generates the unit ideal in $A$. The category of analytic Huber rings is equivalent to the category of analytic Banach rings (over a nonarchimedean field) and it is crucial for the proof to promote analytic Huber rings to analytic Banach rings (and back) to obtain the strictness of the multiplication map by the term defining Laurent coverings. However, $\mathbf{A}_{\textnormal{inf}}(R^{+})$
is never an analytic Huber ring because it always contains a nonempty set of non-analytic valuations corresponding to $\textnormal{Spec}(\mathbf{A}_{\textnormal{inf}}(R^{+})/I)$ where $I$ is the ideal of definition of $\mathbf{A}_{\textnormal{inf}}(R^{+})$. \par

To circumvent the issue of non-analyticity, we show that the streamlined proof of the `stably uniform implies sheafy' argument in \cite{aws} can be formulated slightly more generally for stably uniform Banach rings whose underlying topological ring is Huber:
\begin{Theo}
Let $B$ be a stably uniform Banach ring such that its Banach norm defines a Huber ring $(B, B^{+})$. $\textnormal{Spa}(B, B^{+})$ is an adic space and $\mathcal{O}_{\textnormal{Spa}(B,B^{+})}$ is an acyclic sheaf.
\end{Theo}
In particular, one observes that showing the vanishing of the \v{C}ech cohomologies on arbitrary open coverings on a general affinoid adic space may be first reduced to the case of standard rational coverings and then finally to the case of two term Laurent coverings. This reduction step is in fact true for arbitrary complete Huber rings. $\mathbf{A}_{\textnormal{inf}}(R^{+})$ is a Banach ring with the Gauss norm extended from $(R, R^{+})$. We show that the Gauss norm defines the Huber ring topology on $\mathbf{A}_{\textnormal{inf}}(R^{+})$. Finally to show $\mathbf{A}_{\textnormal{inf}}(R^{+})$ is a stably uniform Banach ring, we use elements from the theory of prismatic cohomology\cite{bs}. We embed $\mathbf{A}_{\textnormal{inf}}(R^{+})$ into its $p$-power-root completion $\mathbf{A}_{\textnormal{inf}}(R^{+})\langle p^{p^{-\infty}}\rangle$ and observe that the latter can be canonically identified as a lens   $\frac{\mathbf{A}_{\textnormal{inf}}(R^{+} [[ T^{p^{-\infty}} ]]) }{([T] - p)}$, i.e. the quotient of the perfect prism $(\mathbf{A}_{\textnormal{inf}}(R^{+}[[T^{p^{-\infty}}]]), ([T]-p))$ by a distinguished element. Furthermore, we can canonically identify rational localizations of $\mathbf{A}_{\textnormal{inf}}(R^{+})\langle p^{p^{-\infty}}\rangle$ with lenses defined by the $\mathbf{A}_{\textnormal{inf}}$ of rational localizations of $R\langle T^{p^{-\infty}}\rangle$:
\begin{Prop}Let $(f_{1},...,f_{n},g)$ be an open ideal of $\mathbf{A}_{\textnormal{inf}}(R^{+})$ defining a rational open subset such that $p \nmid f_{1},...,f_{n},g$. 
We have a canonical topological isomorphism \[  \frac{\mathbf{A}_{\textnormal{inf}}(R^{+})\langle p^{p^{-\infty}} \rangle \langle T_1,...,T_n \rangle}{\overline{(f_1 - T_{1}g,...,f_n - T_{n}g)}}\Big[ \frac{1}{g} \Big] \cong \frac{\mathbf{A}_{\textnormal{inf}}(\frac{R^{+} \langle T_1,...,T_n \rangle}{\overline{(f_{10} - T_{1}g_{0},...,f_n - T_{n}g_{0})}}[\frac{1}{g_{0}}] \langle T^{p^{-\infty}}\rangle_{\rho})}{([T] - p)} , \] 
where $\frac{R^{+} \langle T_1,...,T_n \rangle}{\overline{(f_{10} - T_{1}g_{0},...,f_n - T_{n}g_{0})}}[\frac{1}{g_{0}}] \langle T^{p^{-\infty}}\rangle_{\rho}$ is the completion of $\frac{R^{+} \langle T_1,...,T_n \rangle}{\overline{(f_{10} - T_{1}g_{0},...,f_{n0} - T_{n}g_{0})}}[\frac{1}{g_{0}}] [T^{p^{-\infty}}]$ under the weighted Gauss norm $|\sum_{i}a_{i}T^{i}| = \max_{i}\{|a_{i}|\rho^{i} \}$ with $\rho \in (0,1)$.
\end{Prop}
\sloppy The lens isomorphism for rational localizations of $\mathbf{A}_{\textnormal{inf}}(R^{+})$ allows us to extend the power-multiplicative Banach norm on $R^{+}$ to a power-multiplicative Banach norm on the rational localizations of $\mathbf{A}_{\textnormal{inf}}(R^{+})$. Finally, once we have the acyclicity of the structure sheaf of $\textnormal{Spa}(\mathbf{A}_{\textnormal{inf}}(R^{+}),\mathbf{A}_{\textnormal{inf}}(R^{+}))$, we can follow \cite{aws} 1.9 verbatim to show the equivalence of categories between the category of vector bundles over $\textnormal{Spa}(\mathbf{A}_{\textnormal{inf}}(R^{+}),\mathbf{A}_{\textnormal{inf}}(R^{+}))$ and the category of vector bundles over $\textnormal{Spec}(\mathbf{A}_{\textnormal{inf}}(R^{+}))$:
\begin{Theo}
Let $(R, R^{+})$ be an analytic perfectoid ring in char $p$. The functor $$\mathbf{FPMod}_{\mathbf{A}_{\textnormal{inf}}(R^{+})} \rightarrow \mathbf{Vec}_{\textnormal{Spa}(\mathbf{A}_{\textnormal{inf}}(R^{+}), \mathbf{A}_{\textnormal{inf}}(R^{+}))}: \quad M \rightarrow \widetilde{M}$$ is an equivalence of categories, with quasi-inverse $\mathcal{F} \rightarrow \mathcal{F}(\mathbf{A}_{\textnormal{inf}}(R^{+}))$. In particular, every sheaf in $\mathbf{Vec}_{\textnormal{Spa}(\mathbf{A}_{\textnormal{inf}}(R^{+}), \mathbf{A}_{\textnormal{inf}}(R^{+}))}$ is acyclic.
\end{Theo}

\subsection{Overview of the proof}
In chapter 2, we recollect the basic definitions and discuss uniformity for Huber and Banach rings as well as results for analytic rings. \par
In chapter 3, we fix some notations and conventions for $\mathbf{A}_{\textnormal{inf}}(R^{+})$ and its rational localizations for the rest of the paper. We present a few initial reduction steps.\par
In chapter 4, we introduce perfect prisms and lenses from the theory of prismatic cohomology. We show the rational localizations of $\mathbf{A}_{\textnormal{inf}}(R^{+})$ can be identified as lenses and give the explicit topological lens isomorphism.\par
In chapter 5, using the explicit lens isomorphism, we show that the power-multiplicative Banach norm on $(R,R^{+})$ can be extended to a power-multiplicative norm on the rational localizations of $\mathbf{A}_{\textnormal{inf}}(R^{+})$ agreeing with the Huber ring topology. In particular, we show that the weighted Gauss norm defines the correct Huber ring topology on $\mathbf{A}_{\textnormal{inf}}(R^{+})$, i.e. the topology induced from $(R, R^{+})$.\par
In chapter 6, we demonstrate the streamlined proof of `stably uniform implies sheafy' for stably uniform Banach rings whose underlying topological ring is Huber. In particular, we show the sheafiness and acyclicity of $\textnormal{Spa}(\mathbf{A}_{\textnormal{inf}}(R^{+}),\mathbf{A}_{\textnormal{inf}}(R^{+}))$. \par
In chapter 7, we give results on the gluing of vector bundles over $\textnormal{Spa}(\mathbf{A}_{\textnormal{inf}}(R^{+}),\mathbf{A}_{\textnormal{inf}}(R^{+}))$, showing the equivalence of categories of vector bundles over $\textnormal{Spa}(\mathbf{A}_{\textnormal{inf}}(R^{+}),\mathbf{A}_{\textnormal{inf}}(R^{+}))$ and finite projective modules over $\mathbf{A}_{\textnormal{inf}}(R^{+})$.

\section*{Acknowledgment}
The author would like to thank professor Kiran S. Kedlaya for suggesting this problem and providing many very helpful discussions throughout the writing of this paper. The author is deeply grateful for the years of mentorship provided by Professor Kedlaya throughout their undergraduate and PhD studies. The author would also like to thank Zeyu Liu for reading a preliminary draft of this paper and providing helpful feedback. The author is partially supported by NSF grants DMS-1802161 and DMS-2053473 under professor Kedlaya during the preparation of this project. 
\section{Huber rings, Banach rings, analytic rings and uniformity}
In this section we first define Huber rings and Banach rings. We then define the two different notions of uniformity for Huber rings and Banach rings respectively and how they are related in general and under the assumption of analyticity. 

\begin{Dfn} 
A Huber ring is a topological ring $A$ containing an open subring $A_{0}$ carrying the linear topology induced by a finitely generated ideal $I \subset A_{0}$. The ring $A_{0}$ and the ideal $I$ are called the ring of definition and the ideal of definition respectively.\\
A Huber ring $A$ is Tate if it contains a topologically nilpotent unit, i.e. a pseudo-uniformizer.
\end{Dfn}

\begin{Dfn}
A (nonarchimedean commutative) Banach ring is a ring $B$ equipped with a function $|\bullet|: B \rightarrow \mathbb{R}_{\geq 0}$ satisfying the following conditions:
\begin{enumerate}[(a)]
    \item On the additive group of $B$, $|\bullet|$ is a norm (i.e. a nonarchimedean absolute value such that $|x-y| \leq \textnormal{max} \{|x|,|y|\}$ for all $x,y \in B$).
    \item The norm $|\bullet|$ on $B$ is submultiplicative: $|xy| \leq|x||y|$ for all $x,y \in B$
\end{enumerate}
\end{Dfn}
\begin{Dfn}
     For $B$ a Banach ring, the spectral seminorm on $B$ is the function $|\cdot|_{sp} : B \rightarrow \mathbb{R}_{\geq 0}$ given by
     \[
     |x|_{sp} = \lim_{n\rightarrow \infty} |x^{n}|^{1/n}  \quad (x\in B) 
     \]
\end{Dfn}
In the literature sometimes Huber rings and Banach rings may be implicitly assumed to Tate or analytic. In the rest of our paper, we do not make such an assumption and work with the general definitions.
\begin{Dfn}
A Huber ring $A$ is uniform if the subset (subring) of power-bounded elements $A^{\circ}:= \{a \in A |\{a^{n},n=1,2,...\} \ \textnormal{is bounded}\}$ is bounded.
\end{Dfn}
\begin{Dfn}
A Banach ring $B$ is uniform if one of the following equivalent conditions holds:
\begin{enumerate}[(1)]
    \item The norm on $B$ $|\bullet|_{B}$ is equivalent to some power-multiplicative norm.
    \item For any integer $m>1$, there exists $c>0$ such that $|x^{m}| \geq c|x|^{m}$ for all $x \in B$.
    \item There exists $c > 0$ such that $|x|_{sp} \geq c |x|$ for all $x \in B$.

\end{enumerate}
\end{Dfn}
The next lemma shows that uniformity for Banach rings whose underlying topological ring is Huber is a stronger condition than uniformity for Huber rings. 
\begin{Lemma}(\cite{aws} Remark 1.5.13)\\
Let $B$ be a uniform Banach ring such that its underlying topological ring is Huber, i.e. the Banach norm $|\bullet|_{B}$ defines a Huber ring $(B,B^{+})$. Then the following conditions hold ((1) and (2) are in fact equivalent):
\begin{enumerate}[(1)]
    \item The spectral seminorm $|\bullet|_{sp}$ defines the same topology as the norm $|\bullet|_{B}$.
    \item The underlying Huber ring of B is uniform, i.e. $B^{\circ}$ is bounded.
\end{enumerate}
\end{Lemma}

Next we give the definition of a analytic topological ring as defined in \cite{aws}. We warn the reader that this is a completely different notion from the analytic ring in condensed mathematics in \cite{condense}. This term is renamed as "locally Tate" in \cite{kedlayacondense} to avoid the confusion. However, to keep consistency with \cite{aws}, we shall continue to use the term `analytic topological ring' in this paper.
 \begin{Dfn}
A topological ring $A$ is analytic if the set of topologically nilpotent elements $A^{\circ\circ}$ generate the unit ideal of $A$. 
\end{Dfn}
Using the above definition, we may talk about analytic Huber rings and analytic Banach rings. The terminology 'analytic' comes from the notion of an analytic valuation $v$ on a Huber ring $A$ where $\textnormal{Ker}(v)$ does not contain any open ideals of $A$, introduced in \cite{huber2}. A Huber ring $A$ is analytic if and only if all of its valuations are analytic. This is characterized in the following lemma.
\begin{Lemma}(\cite{aws} Lemma 1.1.3)\\
The following conditions on a general Huber ring $(A,A^{+})$ are equivalent:
    \begin{enumerate}[(1)]
        \item The ring $A$ is analytic.
        \item Any ideal of definition in any ring of definition $A_{0}$ generates the unit ideal in $A$.
        \item Every open ideal of $A$ is trivial.
        \item $\textnormal{Spa}(A,A^{+})$ contains no point on whose residue field the induced valuation is trivial.
    \end{enumerate}
\end{Lemma}

As per Remark 1.5.4 in \cite{aws}, one can promote a general Huber ring $A$ to a Banach ring using the norm defined by the ideal of definition $J$ of $A$: $|x| = \inf \{e^{-n} |\ n\in \mathbb{Z} \ \textnormal{such that} \ xJ^{m} \subset J^{m+n} \ \textnormal{for all} \ m \in \mathbb{Z}_{\geq 0}\}$. Conversely, for an analytic Banach ring $B$, analyticity will guarantee the existence of an ideal of definition in a ring of definition, making $B$ a Huber ring. Then one can freely view an analytic Huber ring as an analytic Banach ring and vice versa. Moreover, the notions of uniformity are equivalent for analytic Huber rings and analytic Banach rings:
\begin{Lemma}(\cite{aws} Remark 1.5.13)\\
Let $B$ be an analytic Banach ring such that the Banach norm of $B$ is equivalent to the norm defined by an ideal of definition. Then $B$ is uniform if and only if $B^{\circ}$ is bounded in $B$.
\end{Lemma}
Given the above lemma, the sheafiness theory of Huber rings is done under the analyticity, (or more restrictively Tate), assumption. Because $\mathbf{A}_{\textnormal{inf}}$ is not analytic, we can not work in the above setting. Because the Gauss norm on $\mathbf{A}_{\textnormal{inf}}$ defines the Huber ring topology on $\mathbf{A}_{\textnormal{inf}}$, we will work under the setting of general uniform Banach rings.\par
Finally we define analytic perfectoid rings, over which we shall define  $\mathbf{A}_{\textnormal{inf}}$.
\begin{Dfn}
Let $(A, A^{+})$ be a uniform analytic Huber ring. $(A,A^{+})$ is an analytic perfectoid ring if there exists an ideal of definition $I \subset A^{+}$ such that $p \in I^{p}$ and the Frobenius map $\varphi:A^{+}/I \rightarrow A^{+}/I^{p}$ is surjective.
\end{Dfn}
\section{Notations and initial reductions}
We fix the following notations throughout the rest of the paper.
\begin{Notation} \label{ainf topology}
Let $(R,R^+)$ be a characteristic $p$ analytic perfectoid pair and let $W(R^+)$ be the ring of Witt vectors over $R^+$. Let $x_0,...,x_m\in R^+$ be generators of an ideal of definition in $R^+$ and let $\pi:W(R^+) \rightarrow R^+$ be the natural mod-$p$ projection map. Then $W(R^+)$ is complete and separated for the topology defined by the ideal $I := \pi^{-1}((x_0,...,x_m)) = (p, [\overline{x}_{0}],...,[\overline{x}_{m}])$. $\mathbf{A}_{\textnormal{inf}}(R^{+})$ is defined as $W(R^{+})$ equipped with the $(p, [\overline{x}_{0}],...,[\overline{x}_{m}])$-adic topology. We will assume $W(R^{+})$ is always equipped with the $(p, [\overline{x}_{0}],...,[\overline{x}_{m}])$-adic topology (unless stated otherwise) and will use $W(R^{+})$ and $\mathbf{A}_{\textnormal{inf}}(R^{+})$ interchangeably.
\end{Notation}
\begin{Notation}
Let $U$ be a rational open subspace of  $X:= \textnormal{Spa}(W(R^{+}),W(R^{+}))$. Then $$U =  X \big(\frac{f_1,...,f_n}{g} \big)=\{v\in X|\ v(f_i) \leq v(g) \neq 0,\;\text{for all}\;i\}$$ where $(f_1,...,f_n,g)$ generates an open ideal in $W(R^+)$. We will use $f_{i0}$ and $g_0$ to denote the reduction mod $p$ of the $f_{i}, g$'s. (To alleviate the abundant appearances of overlines later.) 
\end{Notation}
The next two lemmas will allow us to assume $p \nmid f_{1},...,f_{n},g$. This will be important as we will need to consider the rational localizations of $W(R^+)$ modulo $p$.
\begin{Lemma}
For any open ideal $(f_1,...,f_n,g) \subset W(R^+)$ defining the rational open subspace $U$, without the loss of generality, we may assume $g$ is not topologically nilpotent and in particular $g$ is not divisible by $p$. 
\end{Lemma}
\noindent\textit{Proof.} If $g$ is topologically nilpotent, then we have for all $v \in U$, $v(g) \neq 0$. This shows for all $v \in U$, $v$ is a nontrivial valuation and thus $U \subset \textnormal{Spa}(W(R^{+}),W(R^{+}))^{an}$. Then since $\textnormal{Spa}(W(R^{+}),W(R^{+}))^{an}$ is stably uniform, $\mathcal{O}_X(U)$ is uniform and in fact $U$ is stably uniform.
\qedwhite

\begin{Lemma}
For any open ideal $(f_1,...,f_n,g) \subset W(R^+)$ defining the rational open subspace $U$, without the loss of generality, we may assume for $i = 1,...,n$, all of the $f_i$'s are not divisible by $p$.
\end{Lemma} 
\noindent\textit{Proof.} If $f_j$ is divisible by $p$, we can replace $f_j$ by an element that is "very close" to $f_j$ and not divisible by $p$ without changing the rational open subset. This follows from the following lemma.

We note that for an open ideal $(f_1,...,f_n,g) \subset W(R^{+})$, $(f_1,...,f_n,g)$ mod $p$ is open in $R^{+}$.

\begin{Lemma} \label{parameter perturbation}
\begin{enumerate}[(1)]
    \item Let $U =  X \big(\frac{f_1,...,f_n}{g} \big)$ be a rational open subset defined by an open ideal $(f_1,...,f_n,g)$. Then there exists $k > 0$, only depending on the ideal $(f_1,...,f_n,g)$, such that the ideal $(f_1 + [\overline{x}_{0}]^k,f_{2},...,f_n,[\overline{x}_{0}]^k,g)$ is open in $W(R^{+})$ and defines the same rational open subset i.e. $$X \big(\frac{f_1,...,f_n}{g}\big) =  X \big(\frac{f_1 + [\overline{x}_{0}]^k,f_{2},...,f_n, [\overline{x}_{0}]^k}{g} \big).$$
($[\overline{x}_{0}]$ can be any generator of the ideal of definition of $W(R^{+})$ not equal to $p$) 
    \item $(f_1 + [\overline{x}_{0}]^k,f_2,...,f_n,[\overline{x}_{0}]^k,g)$ (mod $p$) $= (f_{10} + \overline{x}_{0}^{k},f_{20},...,f_{n0}, \overline{x}_{0}^{k},g_{0})$ is open in $R^{+}$. 
\end{enumerate}
\end{Lemma} 

\noindent\textit{Proof.} There exists $k$ such  that $I^k \subset (f_1,...,f_n,g)$ by openness and in particular $I^{k} \subset (f_1 + [\overline{x}_{0}]^k,f_2,...,f_n,[\overline{x}_{0}]^k,g)$. Then there exists $a_1,...,a_{n+1}$ such that $[\overline{x}_{0}]^k = a_1 f_1 +...+a_n f_n + a_{n+1}g$. Then for any $v \in U =  X \big(\frac{f_1,...,f_n}{g} \big)$,
\[
v([\overline{x}_{0}]^k) = v(a_1 f_1 +...+a_n f_n + a_{n+1}g) \leq \max (v(a_1 f_1),...,v(a_n f_n), v(a_{n+1}g)) \leq v(g)
\]
\[
v(f_1 + [\overline{x}_{0}]^k) \leq \max (v(f_1), v([\overline{x}_{0}]^k)) \leq v(g)
\]
and for any $v \in X(\frac{f_1 + [\overline{x}_{0}]^k,f_{2},...,f_n, [\overline{x}_{0}]^k}{g} \big)$
\[
v(f_1) \leq \max (v(f_1 + [\overline{x}_{0}]^k), v([\overline{x}_{0}]^k)) \leq v(g)
\]
For (2), recall that $I = (p, [\overline{x}_{0}],...,[\overline{x}_{m}])$. We know $I^{k} \subset (f_1 + [\overline{x}_{0}]^k,f_{2},...,f_n,[\overline{x}_{0}]^k,g)$. This inclusion mod $p$ gives $(\overline{x}_{0},...,\overline{x}_{m})^{k} \subset (f_{10} + \overline{x}_{0}^{k},f_{20},...,f_{n0}, \overline{x}_{0}^{k},g_{0})$ and thus is open in $R^{+}$.   
\qedwhite

Now we embed $W(R^+)$ into its `integral perfectoidization' $W(R^+) \langle p^{p^{-\infty}} \rangle := \widehat{W(R^+)[p^{p^{-\infty}}]}_{I}$. By the next two lemmas, it will suffice to show  $W(R^+) \langle p^{p^{-\infty}} \rangle$ is stably uniform.

\begin{Lemma} \label{perfectoid inclusion}
Let $$W(R^+) \langle p^{p^{-\infty}} \rangle := \widehat{W(R^+)[p^{p^{-\infty}}]}_{I}$$      and $$X^{'} := \textnormal{Spa}(W(R^+)\langle p^{p^{-\infty}} \rangle,W(R^+)\langle p^{p^{-\infty}} \rangle).$$ The natural strict inclusion $W(R^+) \hookrightarrow W(R^+) \langle p^{p^{-\infty}} \rangle$ splits in the category of topological $W(R^+)$-modules and is stable under taking completed tensor products. More specifically,
let $U$ be a rational open subset of $X$ defined by $f_1,...,f_n,g$ and let $U^{'}$ be the rational open subset of $X^{'}$ defined by $f_1,...,f_n,g$, i.e  \[ U =  X \big(\frac{f_1,...,f_n}{g} \big) \quad \mbox{     and     } \quad U^{'} =  X^{'} \big(\frac{f_1,...,f_n}{g} \big) \] Then there exits a unique strict inclusion $\mathcal{O}_{X}(U) \hookrightarrow  \mathcal{O}_{X^{'}}(U^{'})$ which splits in the category of topological $\mathcal{O}_{X}(U)$-modules and is compatible with rational localizations.
\end{Lemma}

\noindent\textit{Proof.} The strictness and splitting of $W(R^+) \hookrightarrow W(R^+) \langle p^{p^{-\infty}} \rangle$ is clear. We know that $\mathcal{O}_{X}(U) = \frac{W(R^{+}) \langle T_1,...,T_n \rangle}{\overline{(f_1 - T_{1}g,...,f_n - T_{n}g)}}[\frac{1}{g}]$ and $\mathcal{O}_{X^{'}}(U^{'}) = \frac{W(R^{+}) \langle p^{p^{-\infty}} \rangle \langle T_1,...,T_n \rangle}{\overline{(f_1 - T_{1}g,...,f_n - T_{n}g)}}[\frac{1}{g}]$. It is clear that the natural map $W(R^+) \hookrightarrow W(R^+) \langle p^{p^{-\infty}} \rangle$, after taking the completed tensor product with $\frac{W(R^{+}) \langle T_1,...,T_n \rangle}{\overline{(f_1 - T_{1}g,...,f_n - T_{n}g)}}[\frac{1}{g}]$,
\[
\frac{W(R^{+}) \langle T_1,...,T_n \rangle}{\overline{(f_1 - T_{1}g,...,f_n - T_{n}g)}}[\frac{1}{g}] \longrightarrow \frac{W(R^{+}) \langle T_1,...,T_n \rangle}{\overline{(f_1 - T_{1}g,...,f_n - T_{n}g)}}[\frac{1}{g}] \widehat{\otimes}_{W(R^{+})} W(R^+)\langle p^{p^{-\infty}} \rangle
\]
splits in the category of topological $\frac{W(R^{+}) \langle T_1,...,T_n \rangle}{\overline{(f_1 - T_{1}g,...,f_n - T_{n}g)}}[\frac{1}{g}]$ -modules and thus is the inclusion map. Since both the source and target have $(p, [\overline{x}_{0}],...,[\overline{x}_{m}])$-adic topology, the inclusion is strict. By the universal property of rational localization, we have an unique isomorphism
\[
\frac{W(R^{+}) \langle T_1,...,T_n \rangle}{\overline{(f_1 - T_{1}g,...,f_n - T_{n}g)}}[\frac{1}{g}] \widehat{\otimes}_{W(R^{+})} W(R^+)\langle p^{p^{-\infty}} \rangle \cong  \frac{W(R^{+})\langle p^{p^{-\infty}} \rangle \langle T_1,...,T_n \rangle}{\overline{(f_1 - T_{1}g,...,f_n - T_{n}g)}}[\frac{1}{g}].
\]
\qedwhite

\section{Perfect prism and lens}
    We show a natural (topological) isomorphism between $\mathcal{O}_{X^{'}}(U^{'})$ and quotient of a perfect prism by a distinguished element, i.e. a lens. Moreover, the later can be identified as the quotient of $\mathbf{A}_{\textnormal{inf}}$ of a rational localization of the char $p$ perfectoid ring $(R\langle T^{p^{-\infty}} \rangle, R^{+}\langle T^{p^{-\infty}} \rangle)$ by a primitive element.
\begin{Dfn}(\cite{K24} Definition 2.1.1)\\
Let a $\delta$-ring be a pair $(A, \delta)$ where $A$ is a commutative ring and $\delta:A\rightarrow A$ is a map of sets that satisfies the following conditions for all $x,y \in A$:
\begin{enumerate}[(a)]
    \item $\delta(1) = 0$
    \item $\delta(xy)=x^{p}\delta(y) + y^{p}\delta(x) + p\delta(x)\delta(y)$
    \item $\delta(x+y) = \delta(x)+\delta(y) + \sum_{i=1}^{p-1}\frac{(p-1)!}{i!(p-i)!}x^{i}y^{p-i}$.
\end{enumerate}
The map $\delta$ is called a $p$-derivation on $A$.
\end{Dfn}

\begin{Lemma}(\cite{K24} Lemma 2.1.3)
\begin{enumerate}
    \item Let $A$ be a commutative ring. Let $\delta:A\rightarrow A$ be a $p$-derivation. Then the map $\phi:A\rightarrow A$ given by
\[
\phi(x) = x^{p} + p\delta(x)
\]
is a ring homomorphism that induces the Frobenius endomorphism on $A/pA$. $\phi$ will be referred to as the associated Frobenius lift on $(A,\delta)$.
\item If $A$ is $p$-torsion-free, then this construction defines a bijection between $p$-derivations on $A$  and Frobenius lifts on $A$.
\end{enumerate}
\end{Lemma}
\begin{Dfn}(\cite{K24} Definition 3.3.1)\\
A $\delta$-ring $(A, \delta)$ is perfect if $\phi$ is an ismorphism on $A/pA$.
\end{Dfn}
\begin{Dfn}(\cite{K24} Definition 5.1.1) \label{distinguished}
Let $A$ be a $\delta$-ring. An element $d$ is distinguished if $(p,d,\delta(d))$ is the unit ideal of $A$.
\end{Dfn}
\begin{Dfn}(\cite{K24} Definition 5.3.1)
 \begin{enumerate}[(1)]
     \item A $\delta$-pair consists of a pair $(A,J)$ in which A is a $\delta$-ring and $J$ is an ideal.
     \item A prism is a $\delta$-pair such that:
     \begin{enumerate}[(a)]
         \item The ideal $J$ defines a Cartier divisor on $\textnormal{Spec}(A)$.
         \item The ring $A$ is derived $(p, J)$-complete.
         \item $p \in J + \phi(J)$
     \end{enumerate}
     \item  A prism $(A,J)$ is perfect if $A$ is a perfect $\delta$-ring. Then $J$ is principal and any generator of $J$ is a distinguished element.(\cite{K24} 7.2.2)
 \end{enumerate}
\end{Dfn}

\begin{Dfn}(\cite{K24} Definition 8.1.1)\\
A lens $S$ is a ring of the form $S = A/J$ for some perfect prism $(A, J)$.
\end{Dfn}
A lens is really a “integral perfectoid" ring,i.e. it is a perfectoid ring without the topological assumptions of being analytic/Tate and uniform in the context of Huber rings. We know that any rational localization of a perfectoid ring is perfectoid. We would like to show the analogous statement “any rational localization of an `integral perfectoid' ring is `integral perfectoid' " for a suitable “integral perfectoid" ring.

Note that $(W(R^{+}[[T^{p^{-\infty}}]]), ([T] - p))$ is a perfect prism and $W(R^+) \langle p^{p^{-\infty}} \rangle = \frac{W(R^{+}[[T^{p^{-\infty}}]])}{([T] - p)}$. Thus $W(R^+) \langle p^{p^{-\infty}} \rangle$ is a lens.
(We remark that it is also true that \[W(R^+) \langle p^{p^{-\infty}} \rangle \cong \frac{W(R^{+} \langle T^{p^{-\infty}} \rangle)}{([T] - p)}\] i.e. we do not need to take the $T$-adic completion as forming the ring of Witt vectors has the effect of taking $p$-adic completion and $T$ and $p$ are identified in the quotient. We need to take the $T$-adic completion for the formalism of prisms and lenses.)
Next we will show for any rational open subset $U^{'} \subset X^{'}$, $\mathcal{O}_{X^{'}}(U^{'})$ is a lens.

\begin{Prop}(\cite{K24} Proposition 8.2.5)\label{checklens1}\\
A commutative ring $S$ is a lens if and only if the following conditions hold.
\begin{enumerate}[(1)]
    \item The ring $S$ is classically $p$-complete and $S/p$ is semiperfect.
    \item The kernel of the map $\theta_{S}:W(S^{\flat}) \rightarrow S$ is principal.
    \item There exists some $\varpi$ such that $\varpi^{p} = pu$ for some unit $u \in S$.
\end{enumerate}
\end{Prop}

\begin{Prop}(\cite{K24} Proposition 8.2.6) \label{checklens2}\\
A $p$-torsion-free commutative ring $S$ is a lens if and only if the following conditions hold.
\begin{enumerate}[(a)]
    \item The ring $S$ is classically $p$-complete and and $S/p$ is semiperfect.
    \item The ring $S$ is $p$-normal: every $x \in S[p^{-1}]$ with $x^{p} \in S$ belongs to $S$.
    \item There exists some $\varpi$ such that $\varpi^{p} = pu$ for some unit $u \in S$.
\end{enumerate}
\end{Prop}

\begin{Lemma} \label{frob iso}
The Frobenius map $$\phi:\frac{W(R^{+})\langle p^{p^{-\infty}} \rangle \langle T_1,...,T_n \rangle}{\overline{(f_1 - T_{1}g,...,f_n - T_{n}g)}}[\frac{1}{g}]/(p^{\frac{1}{p}}) \rightarrow \frac{W(R^{+})\langle p^{p^{-\infty}} \rangle \langle T_1,...,T_n \rangle}{\overline{(f_1 - T_{1}g,...,f_n - T_{n}g)}}[\frac{1}{g}]/(p)$$ is an isomorphism.
\end{Lemma}
\noindent\textit{Proof.} \sloppy By Stacks Project \cite{sp} Tag 0AMS, we have an explicit description of closure of ideals in adic rings. Since $W(R^{+})\langle p^{p^{-\infty}} \rangle \langle T_1,...,T_n \rangle$ has $(p, [\overline{x}_{0}],...,[\overline{x}_{m}])$-adic topology and $R^{+}\langle T^{p^{-\infty}}\rangle \langle T_1,...,T_n\rangle$ has $(x_0,...,x_m)$-adic topology, we have 
\begin{align*}
&\overline{(f_1 - T_{1}g,...,f_n - T_{n}g)} \\
&= \bigcap_{l=0}^{\infty} \Big( (f_1 - T_{1}g,...,f_n - T_{n}g) + (p, [\overline{x}_{0}],...,[\overline{x}_{m}])^{l} \Big)  \subset W(R^{+})\langle p^{p^{-\infty}} \rangle \langle T_1,...,T_n \rangle\\
&\\
&\overline{(f_{10} - T_{1}g_{0},...,f_{n0} - T_{n}g_{0})}\\
&=  \bigcap_{l=0}^{\infty} \Big( (f_{10} - T_{1}g_{0},...,f_{n0} - T_{n}g_{0}) + (x_0,...,x_m)^{l} \Big) \subset R^{+}\langle T^{p^{-\infty}} \rangle \langle T_1,...,T_n \rangle.
\end{align*}

This shows that taking closure of ideals and taking quotient by $p$ are compatible. Then we have 
\begin{align*}
\frac{W(R^{+})\langle p^{p^{-\infty}} \rangle \langle T_1,...,T_n \rangle}{\overline{(f_1 - T_{1}g,...,f_n - T_{n}g)}}[\frac{1}{g}]/(p)
&= \frac{W(R^{+})\langle p^{p^{-\infty}} \rangle \langle T_1,...,T_n \rangle}{\bigcap_{l}((f_1 - T_{1}g,...,f_n - T_{n}g) + (p, [\overline{x}_{0}],...,[\overline{x}_{m}])^{l})}[\frac{1}{g}]/(p) \\
& \cong \frac{R^{+}\langle T^{p^{-\infty}} \rangle \langle T_1,...,T_n \rangle}{\bigcap_{l}((f_{10} - T_{1}g_{0},...,f_{n0} - T_{n}g_{0}) + (x_0,...,x_m)^{l})}[\frac{1}{g_{0}}]/(T) \\
&=\frac{R^{+}\langle T^{p^{-\infty}} \rangle \langle T_1,...,T_n \rangle}{\overline{(f_{10} - T_{1}g_{0},...,f_{n0} - T_{n}g_{0})}}[\frac{1}{g_{0}}]/(T).
\end{align*}
And similarly we have
\[
\frac{W(R^{+})\langle p^{p^{-\infty}} \rangle \langle T_1,...,T_n \rangle}{\overline{(f_1 - T_{1}g,...,f_n - T_{n}g)}}[\frac{1}{g}]/(p^{\frac{1}{p}}) \cong \frac{R^{+}\langle T^{p^{-\infty}} \rangle \langle T_1,...,T_n \rangle}{\overline{(f_{10} - T_{1}g_{0},...,f_{n0} - T_{n}g_{0})}}[\frac{1}{g_{0}}]/(T^{\frac{1}{p}}).
\]
Now $\frac{R^{+}\langle T^{p^{-\infty}} \rangle \langle T_1,...,T_n \rangle}{\overline{(f_{10} - T_{1}g_{0},...,f_n - T_{n}g_{0})}}[\frac{1}{g_{0}}]$ is a perfect ring in char $p$ and it is clear that 
\[
\phi: \frac{R^{+}\langle T^{^{p-\infty}} \rangle \langle T_1,...,T_n \rangle}{\overline{(f_{10} - T_{1}g_{0},...,f_{n0} - T_{n}g_{0})}}[\frac{1}{g_{0}}]/(T^{\frac{1}{p}}) \rightarrow \frac{R^{+}\langle T^{^{p-\infty}} \rangle \langle T_1,...,T_n \rangle}{\overline{(f_{10} - T_{1}g_{0},...,f_{n0} - T_{n}g_{0})}}[\frac{1}{g_{0}}]/(T)
\]
is an isomorphism.
\qedwhite
\\ 

\begin{Prop} \label{islens}
$\mathcal{O}_{X^{'}}(U^{'}) = \frac{W(R^{+})\langle p^{p^{-\infty}} \rangle \langle T_1,...,T_n \rangle}{\overline{(f_1 - T_{1}g,...,f_n - T_{n}g)}}[\frac{1}{g}]$ is a lens.
\end{Prop}
\noindent\textit{Proof.} We first show $\mathcal{O}_{X^{'}}(U^{'})$ is $p$-torsion-free. Because $p$ does not divide $g$, it suffices to show $\frac{W(R^{+})\langle p^{p^{-\infty}} \rangle \langle T_1,...,T_n \rangle}{\overline{(f_1 - T_{1}g,...,f_n - T_{n}g)}}$ is $p$-torsion-free. Since $R^{+}$ is an integral domain, the map $W(R^{+}) \hookrightarrow W(R^{+}[\frac{1}{g_{0}}])$ is injective and $g$ is invertible in $W(R^{+}[\frac{1}{g_{0}}])$. We have an injective ring homomorphism $$\frac{W(R^{+})\langle p^{p^{-\infty}} \rangle \langle T_1,...,T_n \rangle}{\overline{(f_1 - T_{1}g,...,f_n - T_{n}g)}} \hookrightarrow \frac{W(R^{+}[\frac{1}{g_{0}}])\langle p^{p^{-\infty}} \rangle \langle T_1,...,T_n \rangle}{\overline{(f_1 - T_{1}g,...,f_n - T_{n}g)}} \cong W(R^{+}[\frac{1}{g_{0}}])\langle p^{p^{-\infty}} \rangle \langle f_1,...,f_n \rangle.$$ It is clear that $W(R^{+}[\frac{1}{g_{0}}])\langle p^{p^{-\infty}} \rangle \langle f_1,...,f_n \rangle$ is $p$-torsion-free, so  $\mathcal{O}_{X^{'}}(U^{'})$ is $p$-torsion-free.
Next we will check all the conditions of Proposition \ref{checklens2}. (3) is clear. For (1), classical $p$-completeness is clear and semiperfectness follows by Lemma \ref{frob iso}. We only need to check $\mathcal{O}_{X^{'}}(U^{'})$ is $p$-normal. Let $x \in \mathcal{O}_{X^{'}}(U^{'})[p^{-1}]$ with $x^{p} \in \mathcal{O}_{X^{'}}(U^{'})$. Let $k$ be the smallest nonnegative integer such that $p^{\frac{k}{p}}x \in \mathcal{O}_{X^{'}}(U^{'})$. If $k > 0$, we have \[
(p^{\frac{k}{p}}x)^{p} = p^{k}x^{p} \in p^{k}\mathcal{O}_{X^{'}}(U^{'}) \subset p\mathcal{O}_{X^{'}}(U^{'})
\]
By Lemma \ref{frob iso}, the Frobenius map $\phi: \mathcal{O}_{X^{'}}(U^{'})/(p^{\frac{1}{p}}) \rightarrow \mathcal{O}_{X^{'}}(U^{'})/(p)$ is an isomorphism and we have $p^{\frac{k}{p}}x \in p^{\frac{1}{p}}\mathcal{O}_{X^{'}}(U^{'})$. Hence $p^{\frac{k-1}{p}}x \in \mathcal{O}_{X^{'}}(U^{'})$ which is a contradiction. Thus $k = 0$ and $x \in \mathcal{O}_{X^{'}}(U^{'})$. 
\qedwhite
\\

\begin{Lemma} \label{tilt completion}
Let $S$ be a perfect ring in char $p$ and let $f \in S$ be a nonzerodivisor. Then $(S/f)^{\flat} $ is the $f$-adic completion of $S$, i.e.
\[(S/f)^{\flat} = \varprojlim_{\phi} S/f \cong \varprojlim_{n} S/(f^{n}) .\]
\end{Lemma}
\noindent\textit{Proof.} This is clear because $ \{S/(f^{p^m}) \}_{m}$ is cofinal in $\varprojlim_{n} S/(f^{n})$ as $m \rightarrow \infty$.
\qedwhite
\\

\begin{Prop}
We have a canonical topological isomorphism \[\mathcal{O}_{X^{'}}(U^{'}) =  \frac{W(R^{+})\langle p^{p^{-\infty}} \rangle \langle T_1,...,T_n \rangle}{\overline{(f_1 - T_{1}g,...,f_n - T_{n}g)}}\Big[ \frac{1}{g} \Big] \cong \frac{W(\frac{R^{+} \langle T_1,...,T_n \rangle}{\overline{(f_{10} - T_{1}g_{0},...,f_n - T_{n}g_{0})}}[\frac{1}{g_{0}}] \langle T^{p^{-\infty}}\rangle_{\rho})}{([T] - p)} , \] 
where $\frac{R^{+} \langle T_1,...,T_n \rangle}{\overline{(f_{10} - T_{1}g_{0},...,f_n - T_{n}g_{0})}}[\frac{1}{g_{0}}] \langle T^{p^{-\infty}}\rangle_{\rho}$ is the completion of $\frac{R^{+} \langle T_1,...,T_n \rangle}{\overline{(f_{10} - T_{1}g_{0},...,f_{n0} - T_{n}g_{0})}}[\frac{1}{g_{0}}] [T^{p^{-\infty}}]$ under the weighted Gauss norm $|\sum_{i}a_{i}T^{i}| = \max_{i}\{|a_{i}|\rho^{i} \}$ with $\rho \in (0,1)$.
\end{Prop}
\noindent\textit{Proof.} By Proposition \ref{islens}, we know $\mathcal{O}_{X^{'}}(U^{'})$ is a lens. Consider the natural map $\theta: W(\mathcal{O}_{X^{'}}(U^{'})^{\flat}) \rightarrow \mathcal{O}_{X^{'}}(U^{'})$. By the proof of Proposition \ref{checklens1}, $\mathrm{Ker}(\theta) = ([T] - p)$ and $(W(\mathcal{O}_{X^{'}}(U^{'})^{\flat}), ([T] - p))$ is a perfect prism with 
\[ \label{lens iso} \tag{$\dagger$}
\frac{W(\mathcal{O}_{X^{'}}(U^{'})^{\flat})}{([T] - p)} \cong \mathcal{O}_{X^{'}}(U^{'})\]
By \noindent\textit{\cite{K24} Proposition 7.3.3.}, the isomorphism is natural. 

We have the sequence of isomorphisms:
\begin{align}
\mathcal{O}_{X^{'}}(U^{'})^{\flat} 
&= \Big( \frac{W(R^{+})\langle p^{p^{-\infty}} \rangle \langle T_1,...,T_n \rangle}{\overline{(f_1 - T_{1}g,...,f_n - T_{n}g)}}[\frac{1}{g}] \Big) ^{\flat} \nonumber\\
& \cong \Big( \frac{W(R^{+})\langle p^{p^{-\infty}} \rangle \langle T_1,...,T_n \rangle}{\overline{(f_1 - T_{1}g,...,f_n - T_{n}g)}}[\frac{1}{g}]/(p)\Big)^{\flat} \nonumber\\
& \cong \Big( \frac{R^{+}\langle T^{p^{-\infty}} \rangle \langle T_1,...,T_n  \rangle}{\overline{(f_{10} - T_{1}g_{0},...,f_{n0} - T_{n}g_{0})}}[\frac{1}{g_{0}}]/(T) \Big) ^{\flat} \nonumber\\
&= \Big( \frac{R^{+} \langle T_1,...,T_n \rangle}{\overline{(f_{10} - T_{1}g_{0},...,f_{n0} - T_{n}g_{0})}}[\frac{1}{g_{0}}] \langle T^{p^{-\infty}} \rangle / (T)  \Big) ^{\flat} \\
& \cong \frac{R^{+} \langle T_1,...,T_n \rangle}{\overline{(f_{10} - T_{1}g_{0},...,f_n - T_{n}g_{0})}}[\frac{1}{g_{0}}] \langle T^{p^{-\infty}}\rangle_{\rho},  \ \ \rho \in (0,1).
\end{align}
$(1)$ is true because the closure of the quotient ideal $(f_{10} - T_{1}g_{0},...,f_{n0} - T_{n}g_{0})$ is taken in the $(\overline{x}_{0},...,\overline{x}_{m})$-adic topology and all of $f_{10},...,f_{n0},g_{0}$ are not divisible by $T$ as we assumed all of the $f_{i},g$ are not divisible by $p$. Thus the quotient is independent of $T$. \\
By Lemma \ref{tilt completion}, $(1)$ is the $(T)$-adic completion of $\frac{R^{+} \langle T_1,...,T_n \rangle}{\overline{(f_{10} - T_{1}g_{0},...,f_{n0} - T_{n}g_{0})}}[\frac{1}{g_{0}}] \langle T^{p^{-\infty}} \rangle$, which is
$$\reallywidehat{\frac{R^{+} \langle T_1,...,T_n \rangle}{\overline{(f_{10} - T_{1}g_{0},...,f_{n0} - T_{n}g_{0})}}[\frac{1}{g_{0}}] [[ T^{p^{-\infty}}]]}_{(\overline{x}_{0},...,\overline{x}_{m})}$$ 
with the $(T, \overline{x}_{0},...,\overline{x}_{m})$-adic topology. This is really just the completion of $\frac{R^{+} \langle T_1,...,T_n \rangle}{\overline{(f_{10} - T_{1}g_{0},...,f_{n0} - T_{n}g_{0})}}[\frac{1}{g_{0}}] [T^{p^{-\infty}}]$ under the weighted Gauss norm $|\sum_{i}a_{i}T^{i}| = \max_{i}\{|a_{i}|\rho^{i} \}$ with $\rho \in (0,1)$. And thus we get $(2)$.
Now it is clear that the natural isomorphism (\ref{lens iso}) identifies the $(p, [\overline{x}_{0}],...,[\overline{x}_{m}])$-adic topology on both sides.
\qedwhite
\\

\section{Banach norms and the (weighted) Gauss norm on the ring of Witt Vectors}
Next, we will use the lens isomorphism to show $\mathcal{O}_{X^{'}}(U^{'})$ is a uniform Banach ring by extending the power multiplicative norm on $R^{+}$ to  $\mathcal{O}_{X^{'}}(U^{'})^{\flat}$. The power-multiplicative norm on $\mathcal{O}_{X^{'}}(U^{'})^{\flat}$ extends to the Gauss norm and weighted Gauss norm on $W(\mathcal{O}_{X^{'}}(U^{'})^{\flat})$. Finally the weighted Gauss norm on $W(\mathcal{O}_{X^{'}}(U^{'})^{\flat})$ will descend to a power-multiplicative norm on $\mathcal{O}_{X^{'}}(U^{'})$ and defines the $(p, [\overline{x}_{0}],...,[\overline{x}_{m}])$-adic topology on  $\mathcal{O}_{X^{'}}(U^{'})$.  \par

We first show the weighted Gauss norm defines the $(p, [\overline{x}_{0}],...,[\overline{x}_{m}])$-adic topology on the underlying Huber ring. In the following definitions and lemmas, let $(S,S^{+})$ be a Huber pair such that $S$ is a perfect and uniform Banach ring in char $p$ (then it follows $S^{+}$ is also a perfect and uniform Banach ring). 

\begin{Dfn}
Let $W^{b}(S)$ be the subset of $W(S)$ consisting of the series $\sum_{n=0}^{\infty}p^{n}[\overline{x}_{n}]$ for which the set $\{\overline{x}_{n} :n=0,1,...\}$ is bounded in $S$. $W^{b}(S)$ forms a subring of $W(S)$ and $W(S^{+}) \subset W^{b}(S)$. $W^{b}(S)$ is equipped with the topology of uniform convergence in the $Teichm\ddot{u}ller$ coordinates and this topology coincides with the topology induced by the Gauss norm defined below.
\end{Dfn}

\begin{Dfn} 
    \begin{enumerate}[(1)]
        \item The Gauss norm on $W(S^{+})$ and $W^{b}(S)$ is defined by 
\[
\Bigg|\sum_{n=0}^{\infty}p^{n}[\overline{x}_{n}] \Bigg| = \sup \{ |\overline{x}_{n}|:n = 0,1,..., \}.
\]
        \item The weighted Gauss norm on $W(S^{+})$ and $W^{b}(S)$ is defined by 
\[
\Bigg|\sum_{n=0}^{\infty}p^{n}[\overline{x}_{n}] \Bigg|_{\rho} = \sup \{ \rho^{-n}|\overline{x}_{n}|:n = 0,1,..., \}
\]
for some $\rho \in (0,1)$.
    \end{enumerate}
\end{Dfn}

\begin{Lemma}(Argument due to Kedlaya) \label{bounded witt analytic} \\
Let $(R, R^{+})$ be an analytic perfectoid pair in char $p$. Then $W^{b}(R)$ is an analytic Huber ring with a ring of definition $W(R^{+}) \subset W^{b}(R)$ and an ideal of definition $([\overline{x}_{0}],...,[\overline{x}_{m}]) \subset W(R^{+}).$ 
\end{Lemma}
\noindent\textit{Proof.} We will show $([\overline{x}_{0}],...,[\overline{x}_{m}])$ generates the unit ideal in $W^{b}(R)$. Since $(\overline{x}_0,...,\overline{x}_m)$ generate the unit ideal in $R$, there exist $\overline{a}_{0},...,\overline{a}_{m} \in R$ such that 
\[
\overline{a}_{0} \overline{x}_{0} +...+ \overline{a}_{m}\overline{x}_{m} = 1 \ \in R.
\]
Then for $[\overline{a}_{0}],...,[\overline{a}_{m}] \in W^{b}(R)$, we have 
\[ \label{unit equation} \tag{*}
[\overline{a}_{0}\overline{x}_{0}] +...+ [\overline{a}_{m}\overline{x}_{m}] = [\overline{a}_{0}\overline{x}_{0} +... + \overline{a}_{m} \overline{x}_{m}] + \sum_{n = 1}^{\infty} p^{n}[z_{n}] = 1 + \sum_{n = 1}^{\infty} p^{n}[z_{n}] \ \in W^{b}(R)
\]
where, for each $n$, $z_{n}$ is given by a certain universal degree $1$ homogeneous polynomial in variables $y_{0}^{p^{-n}},...,y_{m}^{p^{-n}}$ evaluated at $(\overline{a}_{0}\overline{x}_{0}),...,(\overline{a}_{m}\overline{x}_{m})$. 
Next we leave $1$ on one side of the equation (\ref{unit equation}) and raise everything to the ($m+1$)-th power:
\begin{align*}
1 &= 1^{m+1} \\
&= \Big( [\overline{a}_{0}\overline{x}_{0}] +...+ [\overline{a}_{m}\overline{x}_{m}] - \sum_{n = 1}^{\infty} p^{n}[z_{n}] \Big) ^{m+1}\\
&= \sum_{\substack{i_{0},...,i_{m} \in \mathbb{Z}[p^{p^{-\infty}}]_{\geq 0}\\
                                                  i_{0} +...+i_{m} = m + 1}} p^{n_{i_{0},...,i_{m}}}b_{i_{0},...,i_{m}}[\overline{a}_{0}]^{i_0}...[\overline{a}_{m}]^{i_m}([\overline{x}_{0}]^{i_0}...[\overline{x}_{m}]^{i_m}) \\
&= \sum_{\substack{j_{0},...,j_{m} \in \mathbb{Z}_{\geq 0}\\
                                                 j_{0} +...+j_{m} = m + 1 \ \textnormal{or} \ m}}\\
                                                 &\hspace{0.2in}\Big(\sum_{\substack{i_{0},...,i_{m} \in \mathbb{Z}[p^{p^{-\infty}}]_{\geq 0}\\
                                                  i_{0} +...+i_{m} = m + 1 \\
                                                  \forall k, \lfloor i_{k} \rfloor = j_{k}}} 
                     p^{n_{i_{0},...,i_{m}}}b_{i_{0},...,i_{m}}[\overline{a}_{0}]^{i_0}...[\overline{a}_{m}]^{i_m}[\overline{x}_{0}]^{i_{0} - j_{0}}...[\overline{x}_{m}]^{i_{m} - j_{m}} \Big) [\overline{x}_{0}]^{j_0}...[\overline{x}_{m}]^{j_m}\\
&=  \label{linear combo} \tag{**} \sum_{\substack{j_{0},...,j_{m} \in \mathbb{Z}_{\geq 0}\\
                                                 j_{0} +...+j_{m} = m + 1 \ \textnormal{or} \ m}} c_{j_{0},...,j_{m}}[\overline{x}_{0}]^{j_0}...[\overline{x}_{m}]^{j_m}
\end{align*}
The indices $i_{0},...,i_{m}$ in the above summations are from the family of homogeneous polynomials defining Witt vector addition and the $b_{i_{0},...,i_{m}}$'s are some universal integer coefficients. Because each $c_{j_{0},...,j_{m}}$ becomes a finite sum modulo any powers of $p$, we have $c_{j_{0},...,j_{m}} \in W(R)$. The set $ \{ [\overline{a}_{0}]^{i_0}...[\overline{a}_{m}]^{i_m}: i_{0} +...+i_{m} = m +1, \; i_{0},...,i_{m} \in \mathbb{Z}[p^{p^{-\infty}}]_{\geq 0} \}$ (where $i_{0},...,i_{m}$ are degrees of the homogeneous polynomials defining Witt vector addition) is bounded under the Gauss norm. We observe that each $[\overline{x}_{0}]^{i_{0} - j_{0}}...[\overline{x}_{m}]^{i_{m} - j_{m}}$ has norm less than 1. Thus every $c_{j_{0},...,j_{m}} \in W^{b}(R)$. Then (\ref{linear combo}) writes $1$ as a finite linear combination of the $[\overline{x}_{0}],...,[\overline{x}_{m}]$. Therefore we have 
\[
1 = \sum_{\substack{j_{0},...,j_{m} \in \mathbb{Z}_{\geq 0}\\
                                                 j_{0} +...+j_{m} = m + 1 \ \textnormal{or} \ m}} c_{j_{0},...,j_{m}}[\overline{x}_{0}]^{j_0}...[\overline{x}_{m}]^{j_m} \ \in ([\overline{x}_{0}],...,[\overline{x}_{m}]) 
                                                \ \subset W^{b}(R).
\]
\qedwhite
\\
\begin{Cor}\label{bounded witt analytic}
Let $(R, R^{+})$ be an analytic perfectoid pair in char $p$. Then $W^{b}(R)$ is an analytic Huber ring with a ring of definition $W(R^{+}) \subset W^{b}(R)$ and an ideal of definition $(p, [\overline{x}_{0}],...,[\overline{x}_{m}]) \subset W(R^{+}).$ \\  
\end{Cor} 

\begin{Lemma}(Argument due to Kedlaya) \\
Let $(R, R^{+})$ be an analytic perfectoid pair in char $p$. Then the Gauss norm on $W(R^{+})$ defines the $([\overline{x}_{0}],...,[\overline{x}_{m}])$-adic topology.
\end{Lemma}
\noindent\textit{Proof.}  In one direction, it is clear that if $f \in ([\overline{x}_{0}],...,[\overline{x}_{m}])^{N}$ for $N$ large, the Gauss norm of $f$ will be small. For the other direction, we need to show if $f \in W(R^{+})$ with sufficiently small Gauss norm $|f|$, then $f \in ([\overline{x}_{0}],...,[\overline{x}_{m}])^{N}$ for sufficiently large $N>0$. To do this, it suffices to show:\\
\textbf{(1)} for such $f \in W(R^{+})$ with small Gauss norm, we can write $f$ as a linear combination of $[\overline{x}_{0}],...,[\overline{x}_{m}]$, i.e. $f = a_{0}[\overline{x}_{0}] +...+ a_{n}[\overline{x}_{n}]$ with $a_{i} \in W(R^{+})$.\\
\textbf{(2)} for any $\varepsilon > 0$, there exists $\delta > 0$ such that if $|f| < \delta$, then $\max_{i}\{|a_{i}|\} < \varepsilon$. \\
Given \textbf{(1)} and \textbf{(2)}, \textbf{(2)} shows that if $f$ has very small norm, then each $a_{i}$ also has small norm. Then one can repeat steps \textbf{(1)} and \textbf{(2)} for each $a_{i}$. This proves the other direction. To show \textbf{(1)}, by Lemma \ref{bounded witt analytic}, there exist $b_{0},...,b_{m} \in W^{b}(R)$ such that 
\[ 
b_0[\overline{x}_{0}] +...+ b_{m}[\overline{x}_{m}] = 1 \ \in W^{b}(R).
\]
(Note that the $b_{i}$'s above are fixed once and for all.) Then multiplying the above equation by $f$, we get 
\[ 
f = fb_{0}[\overline{x}_{0}] +...+ fb_{m}[\overline{x}_{m}] \ \in W^{b}(R).
\]
Now if $f \in W(R^{+})$ has sufficiently small Gauss norm,  for $i \in \{0,...,m\}$, we have $fb_{i} \in W(R^{+})$.   
Thus 
\[
f = fb_{0}[\overline{x}_{0}] +...+ fb_{m}[\overline{x}_{m}] \in ([\overline{x}_{0}],...,[\overline{x}_{m}]) \subset W(R^{+}).
\]
This shows \textbf{(1)} with $a_{i} = fb_{i}$. Now for \textbf{(2)}, given any $\varepsilon$, we can take $\delta = \varepsilon / \max_{i}\{|b_{i}|\}$.
\qedwhite

\begin{Cor}
Let $(R, R^{+})$ be an analytic perfectoid pair in char $p$. Then the weighted Gauss norm on $W(R^{+})$ defines the $(p, [\overline{x}_{0}],...,[\overline{x}_{m}])$-adic topology.
\end{Cor}

The next lemma is \cite{kl} Proposition 3.1.7. In \cite{kl} Proposition 3.1.7 the adic Banach $ \mathbb{F}_{p}$-algebras are assumed to contain a topologically nilpotent unit (pseudo-uniformizer) i.e. they are Tate. The statement and proof of the proposition are valid without the pseudo-uniformizer assumption. We relax the hypothesis and reproduce the proof for the sake of completeness.

\begin{Lemma}(\cite{kl} Proposition 3.1.7) \label{char p rational uniform}\\
Let $(A, A^{+})$ be a general Huber pair such that $A$ is a perfect and uniform Banach ring in char $p$ (the Banach norm gives rise to the Huber ring topology). Then any rational localization of $(A, A^{+})$ is also a perfect and uniform Banach ring. 
\end{Lemma}

\noindent\textit{Proof.} Let $(A, A^{+}) \rightarrow (B, B^{+})$ be a rational localization corresponding to a rational open subset $V \subset \textnormal{Spa}(A,A^{+})$. Then $
B = \frac{A \langle T_{1},...,T_{n} \rangle}{\overline{(w_{1} - T_{1}z,...,w_{n} - T_{n}z)}}[z^{-1}] 
$ where $(w_{1},...,w_{n},z) \subset A$ is an open ideal. $B$ is clearly perfect and $B^{+}$ is also perfect because it is integrally closed in $B$. Equip $B^{1/p}$ with the norm $|x|_{B^{1/p}} = |x^{p}|_{B}^{1/p}$. By applying the inverse of the Frobenius map $\phi^{-1}$, raising norms to the $p$-th power, and using that $A$ is perfect and its norm is power-multiplicative (so its norm is unchanged), we have another rational localization $(A, A^{+}) \rightarrow (B^{1/p}, (B^{+})^{1/p})$ representing $V$. Then we have $B^{1/p} \cong \frac{A \langle T^{1/p}_{1},...,T^{1/p}_{n} \rangle}{\overline{(w^{1/p}_{1} - T^{1/p}_{1}z^{1/p},...,w^{1/p}_{n} - T^{1/p}_{n}z^{1/p})}}[z^{-1/p}]$. The inclusion $A\langle T_1,...,T_n\rangle \hookrightarrow A \langle T^{1/p}_{1},...,T^{1/p}_{n} \rangle$ induces a morphism $(B, B^{+}) \rightarrow (B^{1/p}, (B^{+})^{1/p})$ of Huber pairs over $(A,A^{+})$ which must be an isomorphism by the universal property of rational localizations. This shows $B$ and $B^{+}$ are perfect and uniform Banach rings.
\qedwhite
\\

\begin{Prop}
$\frac{R^{+}\langle T_1,...,T_n \rangle   }{\overline{(f_{10} - T_{1}g_{0},...,f_{n0} - T_{n}g_{0})}}[\frac{1}{g_{0}}]\langle T^{p^{-\infty}}\rangle_{\rho}$ with $\rho \in (0,1)$ is a perfect and uniform Banach ring.
\end{Prop}
\noindent\textit{Proof.} Note that $(f_{10},...,f_{n0},g_{0})$ is open in $R^{+}$ (by Lemma \ref{parameter perturbation} (2)). Then $\frac{R^{+}\langle T_1,...,T_n \rangle   }{\overline{(f_{10} - T_{1}g_{0},...,f_{n0} - T_{n}g_{0})}}[\frac{1}{g_{0}}]$ is a rational localization of the Huber pair $(R^{+},R^{+})$. By Lemma \ref{char p rational uniform}, we know that
$\frac{R^{+}\langle T_1,...,T_n \rangle   }{\overline{(f_{10} - T_{1}g_{0},...,f_{n0} - T_{n}g_{0})}}[\frac{1}{g_{0}}]$   is perfect and uniform as a Banach ring. We give $\frac{R^{+}\langle T_1,...,T_n \rangle}{\overline{(f_{10} - T_{1}g_{0},...,f_{n0} - T_{n}g_{0})}}[\frac{1}{g_{0}}]\langle T^{p^{-\infty}}\rangle_{\rho}$ the weighted Gauss norm $|\sum_{i}a_{i}T^{i}| = \max_{i}\{|a_{i}| \rho^{i} \}$. The (weighted) Gauss extension of a power-multiplicative norm is also a power-multiplicative norm by \cite{K13} Lemma 1.7. This shows $\frac{R^{+}\langle T_1,...,T_n \rangle   }{\overline{(f_{10} - T_{1}g_{0},...,f_{n0} - T_{n}g_{0})}}[\frac{1}{g_{0}}]\langle T^{p^{-\infty}}\rangle_{\rho}$ is a uniform Banach ring and perfectness is clear. 
\qedwhite
\\

\begin{Lemma} \label{witt uniform}(\cite{K13} Lemma 4.1 and Corollary 4.2) \\
If $S$ is a uniform Banach ring, then $W(S^{+})$ and $W^{b}(S)$ are uniform Banach rings under both the Gauss norm and the weighted Gauss norm. 
\end{Lemma}

\begin{Cor}
Let $(R, R^{+})$ be an analytic perfectoid pair in char $p$. Let $W^{b}(R)$ be a Huber ring  with the ring of definition $W(R^{+}) \subset W^{b}(R)$ and the ideal of definition $I:=(p,[\overline{x}_{0}],...,[\overline{x}_{m}]) \subset W(R^{+}).$ Then $(W^{b}(R), W(R^{+}))$ is an analytic stably uniform Huber pair. In particular, $\textnormal{Spa}(W^{b}(R), W(R^{+}))$ is a sheafy adic space.
\end{Cor}
\noindent\textit{Proof.} Since $W^{b}(R)$ is an analytic uniform Huber ring by by Lemma \ref{bounded witt analytic} and Lemma \ref{witt uniform},  $\widehat{W^{b}(R)[p^{p^{-\infty}}]_{I}}$ is an analytic perfectoid ring. We have a strict inclusion $W^{b}(R) \hookrightarrow \widehat{W^{b}(R)[p^{p^{-\infty}}]_{I}}$ which splits in the category of topological $W^{b}(R)$-modules and is stable under rational localization. Therefore $W^{b}(R)$ is a stably uniform analytic Huber ring.
\qedwhite

\begin{Dfn}
An element $z = \sum_{n = 0}^{\infty}p^{n}[\overline{z}_{n}] \in W(S^{+})$ is primitive if  $\overline{z}_0$ is topologically nilpotent and $\overline{z}_1$ is a unit in $S^{+}$ i.e. $z = [\overline{z}_{0}] + pz_{1}$ where $z_1$ is a unit in $W(S^{+})$. Note that a primitive element is a distinguished element (Definition \ref{distinguished}) when we give $W(S^{+})$ the structure of a perfect $\delta$-ring.  
\end{Dfn}

\begin{Prop} \label{quotient uniform} (\cite{aws} Corollary 2.6.10) \\
For $z$ primitive, the quotient norm of the Gauss norm on $W^{b}(S)/(z)$ is power-multiplicative. If in addition the norm on $S$ is multiplicative, then the quotient norm on $W^b(S)/(z)$ is multiplicative.
\end{Prop}

Since we are working with $\mathbf{A}_{\textnormal{inf}}$, we really want the $(p, [\overline{y}_{0}],...,[\overline{y}_{m}])$-adic topology on $W(S^{+})$ and thus we really should be considering the weighted Gauss norm on $W(S^{+})$ and $W^{b}(S)$. However, by \cite{aws} Lemma 2.6.9, the quotient norm of the weighted Gauss norm with $\rho \in (0,1)$ and the Gauss norm agree on $W^{b}(S)/(z)$ when $z$ is a primitive element.\\

\begin{Prop}
 $\mathcal{O}_{X^{'}}(U^{'})$ is a uniform Banach ring i.e. there exists a power-multiplicative norm on $\mathcal{O}_{X^{'}}(U^{'})$ that induces the $(p, [\overline{x}_{0}],...,[\overline{x}_{m}])$-adic topology.
\end{Prop}
\noindent\textit{Proof.} $T$ is topologically nilpotent in $\frac{R^{+}\langle T_1,...,T_n \rangle   }{\overline{(f_{10} - T_{1}g_{0},...,f_{n0} - T_{n}g_{0})}}[\frac{1}{g_{0}}]\langle T^{p^{-\infty}}\rangle_{\rho}$: $|T| = \rho < 1$ and we have $|T^{n}| = |T|^{n} = \rho ^{n} \rightarrow 0$ as n $\rightarrow \infty$. Therefore $[T] - p$ is a primitive element in $W \big( \frac{R^{+}\langle T_1,...,T_n \rangle   }{\overline{(f_{10} - T_{1}g_{0},...,f_{n0} - T_{n}g_{0})}}[\frac{1}{g_{0}}]\langle T^{p^{-\infty}}\rangle_{\rho} \big)$. Then the claim follows from Proposition \ref{quotient uniform}.
\qedwhite
\\

\begin{Cor} \label{rational uniform}
$\mathcal{O}_{X}(U)$ is a uniform Banach ring i.e. there exists a power-multiplicative norm on $\mathcal{O}_{X}(U)$ that induces the $(p, [\overline{x}_{0}],...,[\overline{x}_{m}])$-adic topology.
\end{Cor}
\noindent\textit{Proof.} By Lemma \ref{perfectoid inclusion}, the natural inclusion $\mathcal{O}_{X}(U) \hookrightarrow  \mathcal{O}_{X^{'}}(U^{'})$ is strict. This implies the Banach norm on $\mathcal{O}_{X}(U)$ is equivalent to a power-multiplicative norm.
\qedwhite
\\

\begin{Theo}
Let $(R, R^{+})$ be an analytic perfectoid pair in char $p$. $\mathbf{A}_{\textnormal{inf}}(R^{+})$ is a stably uniform Banach ring and thus $(\mathbf{A}_{\textnormal{inf}}(R^{+}), \mathbf{A}_{\textnormal{inf}}(R^{+}))$ a stably uniform Huber pair.
\end{Theo}
\noindent\textit{Proof.} The theorem follows from the corollary above.
\qedwhite

\section{`Stably uniform implies sheafy' for Banach rings whose underlying topological ring is Huber}
We show the `stably uniform implies sheafy' argument in \cite{aws} works for general Banach rings whose underlying topological ring is Huber, without the analyticity or pseudo-uniformizer assumption. In this section we let $B$ be a stably uniform Banach ring whose Banach norm defines a Huber ring $(B, B^{+})$ and show $\textnormal{Spa}(B,B^{+})$ is a sheafy adic space. In particular, we show $\textnormal{Spa}(\mathbf{A}_{\textnormal{inf}}(R^{+}),\mathbf{A}_{\textnormal{inf}}(R^{+}))$ is a sheafy adic space.

\begin{Lemma}(\cite{aws} Lemma 1.6.3) \label{sheafyreduct}\\
Let $C$ be a cofinal family of rational coverings. Let $X = \textnormal{Spa}(A,A^+)$ where $(A, A^+)$ is a Huber pair. Let $\mathcal{F}$ be a presheaf on $X$ with the property that for any open subset $U$, $\mathcal{F}(U)$ is the inverse limit of $\mathcal{F}(V)$ over all rational subspaces $V \subseteq U$. \par
(a) Suppose that for every rational subspace $U$ of $X$ and every covering $\mathfrak{V} \in C(U) $, the
natural map
\[
\mathcal{F}(U) \rightarrow \check{\mathrm{H}}^0(U,\mathcal{F};\mathfrak{V})
\]
is an isomorphism. Then $\mathcal{F}$ is a sheaf. \par

(b) Suppose that $\mathcal{F}$ is sheaf, and that for every rational subspace $U$ of $X$ and every covering $\mathfrak{V} \in C(U)$, we have $\check{\mathrm{H}}^i(U,\mathcal{F};\mathfrak{V}) = 0$ for all $i > 0$. Then $\mathcal{F}$ is acyclic. 
\end{Lemma}

By Lemma \ref{sheafyreduct}, we need to show for any rational open subspace $U$ of $\textnormal{Spa}(W(R^+),W(R^+))$, the \v{C}ech cohomology groups  $\check{\mathrm{H}}^i(U,\mathcal{O};\mathfrak{V})$ vanish for all $i > 0$ and $\check{\mathrm{H}}^0(U,\mathcal{O};\mathfrak{V})$ is isomorphic to $\mathcal{O}(U)$ where $\mathfrak{V}$ belongs to some cofinal family of rational coverings of $U$. \par

To calculate the \v{C}ech cohomologies of $\textnormal{Spa}(\mathbf{A}_{\textnormal{inf}}(R^{+}), \mathbf{A}_{\textnormal{inf}}(R^{+}))$, we can refine any general rational open coverings to standard rational open coverings.

\begin{Lemma}\label{standardcover}
(\cite{huber} Lemma 2.6)\\
Let $(A,A^+)$be a complete Huber pair. Let $(V_j)_{j\in J}$ be an open covering of $\textnormal{Spa}(A,A^+)$. Then there exist $f_1,...,f_n\in A$ generating the unit ideal such that for every $i=0,...,n$, the rational subset $U\Big(\frac{f_1,...,f_n}{f_i}\Big)$is contained in some $V_j$. Such open coverings are called the standard rational coverings defined by $f_1,...,f_n$.
\end{Lemma}

We remark that for the following two lemmas, the analytic assumption on the Huber rings are not needed because we do not need to use arbitrary open ideals in the Huber ring $A$ to generate the desired binary rational open subsets. The wanted binary rational open subsets come from an induction on the number of the parameters $n$ of a standard rational covering $f_{1},...,f_{n}$ generating the unit ideal in $A$ from Lemma \ref{standardcover}.
\begin{Lemma}(\cite{aws} Lemma 1.6.12)\\
For a general Huber pair $(A,A^{+})$, every open covering of a rational subspace of $\textnormal{Spa}(A,A^{+})$ can be refined by some composition of standard binary rational coverings. 
\end{Lemma} 

\begin{Lemma}\label{simplebinary}
(\cite{aws}  Lemma 1.6.13)\\
Let $(A,A^{+})$ be a complete Huber pair. Every open covering of a rational subspace of  $\textnormal{Spa}(A,A^+)$  can be refined by some composition of coverings, each of which is either a simple Laurent covering or a simple balanced covering.
\end{Lemma}

Let $(B, B^{+}) \rightarrow (C, C^{+}) $ be a rational localization. By the above two lemma, we can compute the \v{C}ech cohomologies on a simple Laurent covering or a simple balanced covering. Then we need to show, for every pair $f,g \in C$ with $g \in \{1, 1-f \}$, the following is a (strict) exact sequence:
\[ 
\begin{tikzcd} \label{cechcomplex1} \tag{$\dagger$}
  0 \arrow[r] & C \arrow[r] & \dfrac{C\langle T \rangle}{\overline{(gT-f)}} \bigoplus \dfrac{C\langle T^{-1} \rangle}{\overline{(g-T^{-1}f)}}  \arrow[r] & \dfrac{C\langle T,T^{-1} \rangle}{\overline{(gT-f)}}   \arrow[r] & 0 
\end{tikzcd}
\]
\\

The next lemma is \cite{aws} Lemma 1.5.26. \cite{aws} Lemma 1.5.26 assumes $A$ to be an uniform analytic Huber ring. The analytic assumption is only needed to promote a uniform Huber ring to a uniform Banch ring and is irrelevant to the rest of the proof. For the sake of completeness, we reproduce the proof here.
\begin{Lemma}(\cite{aws} Lemma 1.5.26) \label{strict inclusion}\\
Suppose $A$ is a uniform Banach ring. Choose $x = \sum_{n=0}^\infty x_nT^n \in A \langle T \rangle$ such that the $x_n$'s generate the unit ideal in $A$. Then multiplication by $x$ defines a strict inclusion $A \langle T \rangle \rightarrow  A \langle T \rangle$. In particular, $(x)$ is a closed ideal in $A \langle T \rangle$. (Same for $A \langle T^\pm \rangle $)
\end{Lemma}

\noindent\textit{Proof.} Let $\mathcal{M}(A)$ be the Gelfand Spectrum of $A$. For $\alpha \in \mathcal{M}(A)$ a multiplicative seminorm on $A$, write $\widetilde{\alpha} \in \mathcal{M}(A \langle T \rangle)$ for the Gauss extension. Then $\widetilde{\alpha}$ is the maximal seminorm on $A \langle T \rangle$ restricting to $\alpha$ on $\mathcal{M}(A)$. Since for a general Banach ring $B$, the spectral seminorm of $B$ equals the supremum over $\mathcal{M}(B)$ (\cite{aws}  Lemma 1.5.22), we may compute the spectral seminorm on $A \langle T \rangle$ as the supremum of $\widetilde{\alpha}$ as $\alpha$  runs over $\mathcal{M}(A)$. \par
Choose $x_0,...,x_m$ that generate the unit ideal in $A$, the quantity 
\[
c := \inf_{\alpha \in \mathcal{M}(A)} \{ \max_{i} \{ \alpha(x_0),...,\alpha(x_m)\} \}
\]
is  positive. For all $y \in A\langle T \rangle$, we have
\[
\sup_{\alpha \in \mathcal{M}(A)}\{\widetilde{\alpha}(xy)\} = \sup_{\alpha \in \mathcal{M}(A)}\{\widetilde{\alpha}(x) \widetilde{\alpha}(y)\} \geq c \sup_{\alpha \in \mathcal{M}(A)}\{\widetilde{\alpha}(y)\}
\]
Since $A$ is uniform, the spectral seminorm is a norm on $A$, this shows multiplication by $x$ is a strict inclusion. In general for a strict morphism $f:M \rightarrow N$ of Banach modules, $\textnormal{Im}(f)$ is closed: $\textnormal{Im}(f) \subset N$ is Hausdorff in the subspace topology and thus is also Hausdorff in the quotient topology from $M$. Thus $\textnormal{Ker}(f)$ is closed and $\textnormal{Im}(f)$ must be complete and thus closed as $N$ is Hausdorff. 
\qedwhite

\begin{Prop} \label{ideal closed}
\begin{enumerate}[(1)]
    \item $(gT - f) \subset C\langle T \rangle$ is closed.
    \item $(gT - f) \subset C\langle T, T^{-1} \rangle$ is closed.
    \item  $(gT^{-1} - f) \subset C\langle T^{-1} \rangle$ is closed. 
\end{enumerate}
\end{Prop}
\noindent\textit{Proof.} These statements follow from the previous lemma.
\qedwhite
\\

By Proposition \ref{ideal closed}, the two term \v{C}ech complex (\ref{cechcomplex1}) becomes
\[
\begin{tikzcd} \label{cechcomplex2} \tag{$\dagger \dagger$}
  0 \arrow[r] & C \arrow[r] & \dfrac{ C \langle T \rangle}{(gT-f)} \bigoplus \dfrac{ C \langle T^{-1} \rangle}{(g-T^{-1}f)}  \arrow[r] & \dfrac{ C \langle T,T^{-1} \rangle}{(gT-f)}   \arrow[r] & 0 
\end{tikzcd}
\]
\\
\begin{Prop}
The two term \v{C}ech complex $($\ref{cechcomplex2}$)$ is exact.
\end{Prop}
\noindent\textit{Proof.} We have the following commutative diagram
\[
\begin{tikzcd} 
         & & 0 \arrow[d] & 0 \arrow[d] \\
     & 0 \arrow[d] \arrow[r] & C\langle T \rangle \bigoplus C\langle T^{-1} \rangle \arrow[d, "{\times(gT-f, \ g-T^{-1}f)}"] \arrow[r, "{\bullet +T^{-1} \bullet}"] & C\langle T,T^{-1} \rangle \arrow[d, "{\times(gT-f)}"] \arrow[r] & 0 \\
      0 \arrow[r] & C \arrow[d] \arrow[r]  & C \langle T \rangle \bigoplus C \langle T^{-1} \rangle \arrow[d] \arrow[r, "{\bullet - \bullet}"] & C \langle T,T^{-1} \rangle \arrow[d] \arrow[r] & 0 \\
  0 \arrow[r] &C \arrow[d] \arrow[r]  & \dfrac{C \langle T \rangle}{(gT-f)} \bigoplus \dfrac{C\langle T^{-1} \rangle}{(g-T^{-1}f)} \arrow[d] \arrow[r] & \dfrac{C\langle T,T^{-1} \rangle}{(gT-f)} \arrow[d]  \arrow[r] & 0 \\
  &0 &0 &0
\end{tikzcd}
\]  
in which all the three columns and the first two rows are exact. We can apply the snake lemma to the first two rows to deduce the exactness at the left and middle of the two term complex. The exactness at the right follows by chasing the bottom right square. 
\qedwhite
\\

\begin{Theo}
Let $B$ be a stably uniform Banach ring such that its Banach norm defines a Huber ring $(B, B^{+})$. $\textnormal{Spa}(B, B^{+})$ is an adic space and $\mathcal{O}_{\textnormal{Spa}(B,B^{+})}$ is an acyclic sheaf.
\end{Theo}
\noindent\textit{Proof.} This follows directly from the above proposition.
\qedwhite
\\
\begin{Theo}
Let $(R, R^{+})$ be an analytic perfectoid pair in char $p$. $\textnormal{Spa}(\mathbf{A}_{\textnormal{inf}}(R^{+}), \mathbf{A}_{\textnormal{inf}}(R^{+}))$ is a sheafy adic space and $\mathcal{O}_{\textnormal{Spa}(\mathbf{A}_{\textnormal{inf}}(R^{+}), \mathbf{A}_{\textnormal{inf}}(R^{+}))}$ is an acyclic sheaf.
\end{Theo}
\noindent\textit{Proof.} This follows from the above theorem and the fact that $(\mathbf{A}_{\textnormal{inf}}(R^{+}), \mathbf{A}_{\textnormal{inf}}(R^{+}))$ is a stably uniform Banach ring.
\qedwhite
\\
\begin{Cor}\label{strictexactcech}
The two term \v{C}ech complex $($\ref{cechcomplex2}$)$ is strict exact.
\end{Cor}
\noindent\textit{Proof.} The first map is an isometry for the spectral seminorm by \cite{aws}  Remark 1.5.25 and thus is a strict inclusion by uniformity. For the second map we consider the bottom right corner of the commutative diagram in the last proposition.
\[
\begin{tikzcd}
    C\langle T \rangle \bigoplus C\langle T^{-1} \rangle \arrow[d, two heads] \arrow[r, two heads, "{\bullet - \bullet}"] & C\langle T,T^{-1} \rangle \arrow[d, two heads] \\
    \dfrac{C\langle T \rangle}{(gT-f)} \bigoplus \dfrac{C\langle T^{-1} \rangle}{(g-T^{-1}f)} \arrow[r] & \dfrac{C\langle T,T^{-1} \rangle}{(gT-f)}  
\end{tikzcd}
\]
All four modules (rings) in the above diagram are complete in the $(p, [\overline{x}_{0}],...,[\overline{x}_{m}])$-adic topology. It is clear that the top horizontal map and the two vertical surjections map open basis to open subsets of the targets. Thus the top horizontal map and the two vertical surjections are open (the ideals in the quotients are closed). Thus the bottom horizontal map is open and strict. Since both maps in the short exact sequence are strict, the strictness of the middle part is clear. 
\qedwhite 

\section{Gluing finite projective modules over $\mathbf{A}_{\textnormal{inf}}(R^{+})$}
Since the sheafiness of $\textnormal{Spa}(\mathbf{A}_{\textnormal{inf}}(R^{+}), \mathbf{A}_{\textnormal{inf}}(R^{+}))$ is established, we follow \cite{aws} 1.9 to show the equivalence of categories between the category of vector bundles over $\textnormal{Spa}(\mathbf{A}_{\textnormal{inf}}(R^{+}), \mathbf{A}_{\textnormal{inf}}(R^{+}))$ and the category of finite projective modules over $\mathbf{A}_{\textnormal{inf}}(R^{+})$. Recently, using condensed mathematics, Kedlaya showed the equivalence of categories between the category of vector bundles over $\textnormal{Spa}(A,A^{+})$ and and the category of finite projective modules over $A$ for a general sheafy Huber pair $(A, A^{+})$\cite{kedlayacondense}. Here we give a direct proof.\par

\begin{Dfn}
Let $(A, A^{+})$ be a Huber pair. For any $A$-module $M$, let $\widetilde{M}$ be the presheaf on $\textnormal{Spa}(A,A^{+})$ such that for open $V \subset \textnormal{Spa}(A,A^{+})$, 
\[
\widetilde{M}(V) = \varprojlim_{\textnormal{Spa}(B,B^{+}) \subset V }M \widehat{\otimes} B,
\]
where the inverse limit is taken over all rational localizations $(A,A^{+}) \rightarrow (B,B^{+})$
\end{Dfn}

\begin{Dfn}
Let $(A,A^{+})$ be a Huber pair. Let $\mathbf{FPMod}_{A}$ denote the category of finite projective $A$-modules. A vector bundle on $\textnormal{Spa}(A,A^{+})$ is a sheaf  $\mathcal{F}$ of $\mathcal{O}_{\textnormal{Spa}(A,A^{+})}$-modules which is locally of the form $\widetilde{M}$ for finite projective $A$-module $M$. More specifically, there exists a finite covering $\big\{ V_{i} \big\}_{i=1}^{n}$ of $\textnormal{Spa}(A,A^{+})$ by rational subspaces such that for each $i$, $M_{i} := \mathcal{F}(V_{i}) \in \mathbf{FPMod}_{\mathcal{O}(V_{i})}$ and the canonical morphism $\widetilde{M}_{i} \rightarrow \mathcal{F}|_{V_{i}}$ of sheaves of $\mathcal{O}|_{V_{i}}$-modules is an isomorphism.
Let $\mathbf{Vec}_{\textnormal{Spa}(A,A^{+})}$ denote the category of vector bundles on $\textnormal{Spa}(A,A^{+})$.The functor $\mathbf{FPMod}_{A} \rightarrow \mathbf{Vec}_{\textnormal{Spa}(A,A^{+})}: M \rightarrow \widetilde{M}$ is exact by the flatness of finite projective modules.
\end{Dfn}

\begin{Theo} \label{acyclicmodule}
For any finite projective $\mathbf{A}_{\textnormal{inf}}(R^{+})$-modules $M$, the presheaf $\widetilde{M}$ is an acyclic sheaf.
\end{Theo}
\noindent\textit{Proof.} Since $M$ is a direct summand of a finite free $\mathbf{A}_{\textnormal{inf}}(R^{+})$-module, we may reduce to the case $M=\mathbf{A}_{\textnormal{inf}}(R^{+})$. By results in section 4, we may check the claim on binary standard rational coverings. Then the theorem follows from Proposition 4.7.
\qedwhite 
\par
Following \cite{aws} Remark 1.6.16, given that $\mathbf{A}_{\textnormal{inf}}(R^{+})$ is sheafy and $\widetilde{M}$ is acyclic for every finite projective $\mathbf{A}_{\textnormal{inf}}(R^{+})$-module, it suffices to consider a bundle which is specified by modules on each term of a composition of simple Laurent coverings and simple balanced coverings by Lemma \ref{simplebinary}. \\
We continue to fix the notation that $U$ is a rational open subspace of $\textnormal{Spa}(\mathbf{A}_{\textnormal{inf}}(R^{+}), \mathbf{A}_{\textnormal{inf}}(R^{+}))$ with ring of global sections 
$\mathcal{O}_X(U)$ and consider the two term coverings $U\big(\frac{f}{g}\big)$ and $U\big(\frac{g}{f}\big)$ for every pair $f, g \in \mathcal{O}_X(U)$ with $g \in \{1, 1-f\}$. Finally we let $\mathcal{O}_X(U)\langle \frac{f}{g}\rangle := \dfrac{ \mathcal{O}_{X}(U)\langle T \rangle}{(gT-f)}$, $\mathcal{O}_X(U)\langle \frac{g}{f}\rangle := \dfrac{ \mathcal{O}_{X}(U)\langle T^{-1} \rangle}{(g-T^{-1}f)}$, and $\mathcal{O}_X(U)\langle \frac{f}{g}, \frac{g}{f}\rangle := \dfrac{ \mathcal{O}_{X}(U)\langle T,T^{-1} \rangle}{(gT-f)}$.\\
\begin{Lemma}\label{matrixfactorize}
(\cite{kl} Lemma 2.7.2) \\
Let $R_{1} \rightarrow S$, $R_{2} \rightarrow S$ be bounded homomorphisms of Banach rings (not necessarily containing topologically nilpotent units) such that the sum homomorphism $\phi: R_{1} \oplus R_{2} \rightarrow S$ of groups is strict surjective. Then there exists a constant $c > 0$ such that for every positive integer $n$, every matrix $U \in GL_{n}(S)$ with $|U - 1| < c$ can be written in the form $\phi(U_{1})\phi(U_{2})$ with $U_{i} \in GL_{n}(R_{i})$, $i=1,2$.
\end{Lemma}
The next lemma is \cite{aws} 1.9.6, where finite projective modules and pseudo-coherent modules over analytic Huber rings are treated in parallel. We reproduce the proof to highlight the independence from analyticity and the open mapping theorem.
\begin{Lemma}\label{fiberproduct}
(\cite{aws} 1.9.6)\\
Let $M_{1}$,$M_{2}$,$M_{12}$ be finitely generated modules over $\mathcal{O}_X(U)\langle \frac{f}{g}\rangle$, $\mathcal{O}_X(U)\langle \frac{g}{f}\rangle $, $\mathcal{O}_X(U)\langle\frac{f}{g},\frac{g}{f}\rangle$ respectively. Let $$\psi_{1}: M_{1} \otimes_{\mathcal{O}_X(U)\langle \frac{f}{g}\rangle} \mathcal{O}_X(U)\langle\frac{f}{g},\frac{g}{f}\rangle \rightarrow M_{12}, \qquad \psi_{2}: M_{2} \otimes_{\mathcal{O}_X(U)\langle \frac{g}{f}\rangle} \mathcal{O}_X(U)\langle\frac{f}{g},\frac{g}{f}\rangle \rightarrow M_{12}$$ be isomorphisms. Then we have \\
$\quad$(a) The map $\psi : M_{1} \oplus M_{2} \rightarrow M_{12}$ taking $(\textbf{v}, \textbf{w})$ to $\psi_{1}(\textbf{v}) - \psi_{2}(\textbf{w})$ is strict surjective. \\
$\quad$(b) For $M := ker(\psi)$, the induced maps \[
M \otimes_{\mathcal{O}_X(U)} \mathcal{O}_X(U)\langle \frac{f}{g}\rangle \rightarrow M_{1}, \qquad M \otimes_{\mathcal{O}_X(U)} \mathcal{O}_X(U)\langle \frac{g}{f}\rangle \rightarrow M_{2}
\]
$\; \qquad$are strict surjective.
\end{Lemma}
\noindent\textit{Proof.} Let $\textbf{v}_{1},...,\textbf{v}_{n}$ and $\textbf{w}_{1},...,\textbf{w}_{n}$ be generating sets of $M_{1}$ and $M_{2}$ respectively of the same cardinality. We chose $n \times n$ matrices $V$ and $W$ over $\mathcal{O}_X(U)\langle\frac{f}{g},\frac{g}{f}\rangle$ such that 
$\psi_{2}(\textbf{w}_{j}) = \sum_{i}V_{ij}\psi_{1}(\textbf{v}_{i})$ and $\psi_{1}(\textbf{v}_{j}) = \sum_{i}W_{ij}\psi_{2}(\textbf{w}_{i})$.\par
By Corollary \ref{strictexactcech}, the map $\phi:\mathcal{O}_X(U)\langle \frac{f}{g}\rangle \oplus \mathcal{O}_X(U)\langle \frac{g}{f}\rangle \rightarrow \mathcal{O}_X(U)\langle\frac{f}{g},\frac{g}{f}\rangle$ is strict surjective and there exists $c > 0$ such that Lemma \ref{matrixfactorize} holds for our $\phi$. Since $\mathcal{O}_X(U)\langle \frac{f}{g}\rangle[f^{-1}]$ is dense in $\mathcal{O}_X(U)\langle\frac{f}{g},\frac{g}{f}\rangle$, we can choose a nonnegative integer $m$ and an $n \times n$ matrix $W^{'}$ over $\mathcal{O}_X(U)\langle \frac{g}{f}\rangle$ so that $|V(f^{-m}W^{'} - W)| < c$. Then we have $1 + V(f^{-m}W^{'} - W) = \phi(X_{1})\phi(X_{2}^{-1})$ with $X_{1} \in GL_{n}(\mathcal{O}_X(U)\langle \frac{f}{g}\rangle), X_{2} \in GL_{n}(\mathcal{O}_X(U)\langle \frac{g}{f}\rangle)$. We shall omit $\phi$ and write  $1 + V(f^{-m}W^{'} - W) = X_{1} X_{2}^{-1}$ for clarity.\par
Now define elements $\textbf{x}_{j} \in M_{1} \oplus M_{2}$ by the formula \[
\textbf{x}_{j} = (\textbf{x}_{j,1}, \textbf{x}_{j,2}) = \Big( \sum_{i}f^{m}(X_{1})_{ij}\textbf{v}_{i}, \sum_{i}(W^{'}X_{2})_{ij}\textbf{w}_{i}) \Big), \qquad j = 1,...,n.
\]
Then for every $j$,\[
\psi_{1}(\textbf{x}_{j,1}) - \psi_{2}(\textbf{x}_{j,2}) = \sum_{i}(f^{m}X_{1} - VW^{'}X_{2})_{ij}\psi_{1}(\textbf{v}_{i}) = \sum_{i}f^{m}((1-VW)X_{2})_{ij}\psi_{1}(\textbf{v}_{i}) = 0,
\]
since $VW$ is the identity and therefore $\textbf{x}_{j} \in M$. Since $X_1 \in GL_{n}(\mathcal{O}_X(U)\langle \frac{f}{g}\rangle)$ and $\{\textbf{v}_{1},...,\textbf{v}_{n}\}$ is a generating set of $M_{1}$, we see that the map $M \otimes_{\mathcal{O}_X(U)}\mathcal{O}_X(U)\langle \frac{f}{g}\rangle \rightarrow M_{1}$ induces a strict surjection onto $f^{m}M_{1}$. \par
The induced map $M \otimes_{\mathcal{O}_X(U)}\mathcal{O}_X(U)\langle \frac{f}{g}, \frac{g}{f}\rangle \rightarrow M_{12}$ is strict surjective as $f$ is invertible in $\mathcal{O}_X(U)\langle \frac{f}{g}, \frac{g}{f}\rangle$. By Corollary 4.9 and tensoring $M$ with the exact two term \v{C}ech complex $($\ref{cechcomplex2}$)$, we get a strict surjection $M \otimes_{\mathcal{O}_X(U)}(\mathcal{O}_X(U)\langle \frac{f}{g}\rangle \oplus \mathcal{O}_X(U)\langle \frac{g}{f}\rangle) \rightarrow M_{12}$ which factors through $\psi$. Thus we obtain (a). \par
For each $\textbf{v} \in M_{2}$, $\psi_{2}(\textbf{v})$ lifts to $M \otimes_{\mathcal{O}_X(U)}(\mathcal{O}_X(U)\langle \frac{f}{g}\rangle \oplus \mathcal{O}_X(U)\langle \frac{g}{f}\rangle)$ by above, then we can find $\textbf{w}_{1} \in M_{1}$, $\textbf{w}_{2} \in M_{2}$ in the images of the base extension maps from $M$ with $\psi_{1}(\textbf{w}_{1}) - \psi_{2}(\textbf{w}_{2}) = \psi_{2}(\textbf{v})$. Then $(\textbf{w}_{1}, \textbf{v} + \textbf{w}_{2}) \in M$, then $\textbf{w}_{2}, \textbf{w}_{2} + \textbf{v} \in \textnormal{Im}(M \otimes_{\mathcal{O}_X(U)} \mathcal{O}_X(U)\langle \frac{g}{f}\rangle \rightarrow M_{2})$ and therefore $M \otimes_{\mathcal{O}_X(U)} \mathcal{O}_X(U)\langle \frac{g}{f}\rangle \rightarrow M_{2}$ is stric surjective. Swapping the role of $f, g$, we get that $M \otimes_{\mathcal{O}_X(U)} \mathcal{O}_X(U)\langle \frac{f}{g}\rangle \rightarrow M_{1}$ is strict surjective. This show (b).
\qedwhite 
\\
\begin{Lemma}\label{maxidealker}
(\cite{wedhorn} Lemma 7.51)\\
Let $(A,A^{+})$ be a complete Huber ring and let $\mathfrak{m} \subset A$ be a maximal ideal. Then $\mathfrak{m}$ is closed and there exists $v \in \textnormal{Spa}(A, A^{+})$ with $\textrm{Ker}(v) = \mathfrak{m}$.
\end{Lemma}
\noindent\textit{Proof.} $A^{\circ \circ}$, the ideal of the set of all topological nilpotent elements in $A$, is open and thus $1 + A^{\circ \circ}$ is open. $1 + A^{\circ \circ} \subset A^{\times}$ implies $A^{\times}$ is open and thus $A \setminus A^{\times}$ is closed. Therefore $\mathfrak{m}$ is closed and $A/ \mathfrak{m}$ is Hausdorff which implies $\textnormal{Spa}(A/ \mathfrak{m}) \neq \emptyset$.
\qedwhite
\begin{Lemma}(\cite{aws} Lemma 1.9.8)\label{maxidealsurjection}\\
The image of the natural map $\textnormal{Spec}(\mathcal{O}_X(U)\langle \frac{f}{g}\rangle \oplus \mathcal{O}_X(U)\langle \frac{g}{f}\rangle) \rightarrow \textnormal{Spec}(\mathcal{O}_X(U))$ contains all maximal ideals of $\mathcal{O}_X(U)$.
\end{Lemma}
\noindent\textit{Proof.} By Lemma \ref{maxidealker}, for every maximal ideal $\mathfrak{m} \subset \mathcal{O}_X(U)$, there is $v \in \textnormal{Spa}(\mathcal{O}_X(U))$ with $\textnormal{Ker}(v) = \mathfrak{m}$. Since $f,g$ forms a binary covering on $\textnormal{Spa}(\mathcal{O}_X(U))$, $v$ extends/factorizes through one of $\mathcal{O}_X(U)\langle \frac{f}{g}\rangle$ or $\mathcal{O}_X(U)\langle \frac{g}{f}\rangle$. The kernel of the extension of $v$ is a prime ideal contracting to $\mathfrak{m}$.
\qedwhite
\\
\begin{Prop} \label{essentialsurjective} Let $U$ be a rational open subspace of $\textnormal{Spa}(\mathbf{A}_{\textnormal{inf}}(R^{+}), \mathbf{A}_{\textnormal{inf}}(R^{+}))$.
There is an exact functor \[
    \mathbf{FPMod}_{\mathcal{O}_X(U)\langle \frac{f}{g}\rangle} \times_{\mathbf{FPMod}_{\mathcal{O}_X(U)\langle\frac{f}{g},\frac{g}{f}\rangle}} \mathbf{FPMod}_{\mathcal{O}_X(U)\langle \frac{g}{f}\rangle} \rightarrow \mathbf{FPMod}_{\mathcal{O}_X(U)}
\]
given by taking equalizers. The composition of this functor with the base extension functor in the opposite direction is naturally isomorphic to the identity.
\end{Prop}
\noindent\textit{Proof.} Using the notations from Lemma \ref{fiberproduct}, if $M_{1}$,$M_{2}$,$M_{12}$ are finite projective modules over $\mathcal{O}_X(U)\langle \frac{f}{g}\rangle$, $\mathcal{O}_X(U)\langle \frac{g}{f}\rangle $, $\mathcal{O}_X(U)\langle\frac{f}{g},\frac{g}{f}\rangle$ respectively, then we need to show $M$ is finite projective over $\mathcal{O}_X(U)$ and the maps \begin{equation}
\label{projectivebasechange}
M \otimes_{\mathcal{O}_X(U)} \mathcal{O}_X(U)\langle \frac{f}{g}\rangle \rightarrow M_{1}, \qquad M \otimes_{\mathcal{O}_X(U)} \mathcal{O}_X(U)\langle \frac{g}{f}\rangle \rightarrow M_{2}
\end{equation} 

are isomorphisms. By Lemma \ref{fiberproduct}, we can choose a finite free $\mathcal{O}_X(U)$-module $F$ and a (not necessarily surjective) $\mathcal{O}_X(U)$-linear map $F \rightarrow M$ such that for $F_{1}, F_{2}, F_{12}$ the respective base extensions of $F$, the induced maps \[
F_{1} \rightarrow M_{1}, \qquad F_{2} \rightarrow M_{2}, \qquad F_{12} \rightarrow M_{12}
\]
are surjective. Let \[
N_{1}:= \textnormal{Ker}(F_{1} \rightarrow M_{1}), \quad 
N_{2}:= \textnormal{Ker}(F_{2} \rightarrow M_{2}), \quad
N_{12}:= \textnormal{Ker}(F_{12} \rightarrow M_{12})
\]
and put $N = \textnormal{Ker}(N_{1} \oplus N_{2} \rightarrow N_{12})$. Consider the following commutative diagram with the second and third columns exact:
\[
\begin{tikzcd} 
& 0 \arrow[d] & 0 \arrow[d] & 0 \arrow[d] \\
 0 \arrow[r] & N \arrow[d] \arrow[r] & N_{1} \oplus N_{2} \arrow[d] \arrow[r] & N_{12} \arrow[d] \arrow[r, dashed] & 0 \\
0 \arrow[r] & F \arrow[d] \arrow[r]  & F_{1} \oplus F_{2} \arrow[d] \arrow[r] & F_{12} \arrow[d] \arrow[r] & 0 \\
  0 \arrow[r] & M \arrow[d, dashed] \arrow[r]  & M_{1} \oplus M_{2} \arrow[d] \arrow[r] & M_{12} \arrow[d]  \arrow[r] & 0 \\
  &0 &0 &0 
\end{tikzcd}
\]  
The exactness of the first column (minus the dashed arrows) is obtained by applying the snake lemma to the second and third columns (injectivity of $N \rightarrow F$ follows by chasing the top left square). We have $N= \textnormal{Ker}(F \rightarrow M)$. \par
Because $N_{1}, N_{2}, N_{12}$ are finite projective modules, we have isomorphisms \[
 N_{1} \otimes_{\mathcal{O}_X(U)\langle \frac{f}{g}\rangle} \mathcal{O}_X(U)\langle\frac{f}{g},\frac{g}{f}\rangle \cong N_{12}, \qquad N_{2} \otimes_{\mathcal{O}_X(U)\langle \frac{g}{f}\rangle} \mathcal{O}_X(U)\langle\frac{f}{g},\frac{g}{f}\rangle \cong N_{12}.
\]
Therefore the modules $N_{1}, N_{2}, N_{12}$ form an object of the fiber product category. By Lemma \ref{fiberproduct} again, we see the top right dash arrow $N_{1} \oplus N_{2} \rightarrow N_{12}$ is exact and applying the snake lemma to the second and third columns of the diagram again gives the exactness of the bottom left dash arrow $F \rightarrow M$.\par
Consider the diagram
\[
\begin{tikzcd} 
& N \otimes_{\mathcal{O}_X(U)} \mathcal{O}_X(U)\langle \frac{f}{g}\rangle  \arrow[d] \arrow[r] & F \otimes_{\mathcal{O}_X(U)} \mathcal{O}_X(U)\langle \frac{f}{g}\rangle \arrow[d] \arrow[r] & M \otimes_{\mathcal{O}_X(U)} \mathcal{O}_X(U)\langle \frac{f}{g}\rangle  \arrow[d] \arrow[r] & 0 \\
 0 \arrow[r] & N_{1} \arrow[r] & F_{1} \arrow[r] & M_{1} \arrow[r] & 0 \\
\end{tikzcd}
\]  
with exact rows. By Lemma \ref{fiberproduct}, both outside vertical arrows are surjective. Since the middle arrow is an isomorphism, the five lemma gives the injectivity of the right vertical arrow, i.e., $M \otimes_{\mathcal{O}_X(U)} \mathcal{O}_X(U)\langle \frac{f}{g}\rangle \cong M_{1}$. Replacing $M_{1}$ with $M_{2}$, we see the maps in (\ref{projectivebasechange}) are isomorphisms. \par
Now it remains to show $M \in \mathbf{FPMod_{ \mathcal{O}_X(U)}}$. We have shown $M$ is a finitely presented $ \mathcal{O}_X(U)$-module and $M_{\mathfrak{m}}$ is a finite free $ \mathcal{O}_X(U)_{\mathfrak{m}}$-module for every maximal ideal $\mathfrak{m}$ of $ \mathcal{O}_X(U)$ by Lemma \ref{maxidealsurjection} and the isomorphisms in (\ref{projectivebasechange}). By Stacks Project\cite{sp} Tag 00NX, $M$ is a finite projective $ \mathcal{O}_X(U)$-module.
\qedwhite
\\
\begin{Theo}
Let $(R, R^{+})$ be an analytic perfectoid pair in char $p$. The functor $$\mathbf{FPMod}_{\mathbf{A}_{\textnormal{inf}}(R^{+})} \rightarrow \mathbf{Vec}_{\textnormal{Spa}(\mathbf{A}_{\textnormal{inf}}(R^{+}), \mathbf{A}_{\textnormal{inf}}(R^{+}))}: \quad M \rightarrow \widetilde{M}$$ is an equivalence of categories, with quasi-inverse $\mathcal{F} \rightarrow \mathcal{F}(\mathbf{A}_{\textnormal{inf}}(R^{+}))$. In particular, by Theorem \ref{acyclicmodule}, every sheaf in $\mathbf{Vec}_{\textnormal{Spa}(\mathbf{A}_{\textnormal{inf}}(R^{+}), \mathbf{A}_{\textnormal{inf}}(R^{+}))}$ is acyclic.
\end{Theo}
\noindent\textit{Proof.} It suffices to show for any rational subspace $U \subset \textnormal{Spa} (\mathbf{A}_{\textnormal{inf}}(R^{+}), \mathbf{A}_{\textnormal{inf}}(R^{+}))$ and any binary standard coverings of $U$ by $U\big(\frac{f}{g}\big)$ and $U\big(\frac{g}{f}\big)$ for  $f, g \in \mathcal{O}_X(U)$ with $g \in \{1, 1-f\}$,
\[ \mathbf{FPMod}_{\mathcal{O}_X(U)} \rightarrow 
    \mathbf{FPMod}_{\mathcal{O}_X(U)\langle \frac{f}{g}\rangle} \times_{\mathbf{FPMod}_{\mathcal{O}_X(U)\langle\frac{f}{g},\frac{g}{f}\rangle}} \mathbf{FPMod}_{\mathcal{O}_X(U)\langle \frac{g}{f}\rangle} 
\]
is an exact equivalence of categories. The functor is fully faithful by Theorem \ref{acyclicmodule}, exact by flatness of finite projective modules, and essentially surjective by Proposition \ref{essentialsurjective}.
\qedwhite

\bibliography{ref}
\bibliographystyle{alpha}

\end{document}